\documentclass[11pt]{article}
\usepackage{geometry}                
\usepackage{graphicx}
\usepackage{amssymb}
\usepackage{epstopdf}
\usepackage[applemac]{inputenc}
\usepackage{latexsym,bm}
\usepackage{mathrsfs,amsmath,amssymb,amsthm}
\usepackage{enumitem}

 \usepackage{graphicx}
\newtheorem{theo}{Theorem}[section]
\newtheorem{defi}{Definition}[section]
\newtheorem{lem}{Lemma}[section]
\newtheorem{propo}{Proposition}[section]
\newtheorem{propri}{Property}[section]
\newtheorem{coro}{Corollary}[section]
\newtheorem{rem}{Remark}[section]

\newcommand{\modint}{\displaystyle\copy\tratto\kern-10.4pt\int\limits}

\def\RR{\mathcal R}
\def\mm{\mu}
\def\M{\mathcal M}

\def\R{{\rm I\!R}}
\def\N{{\rm I\!N}}

\def\OV{\overline}

\def\G{\Gamma}
\def\O{\Omega}
\def\log{{\rm Log\,}}
\def\a{\alpha}
\def\b{\beta}

\def\eps{\varepsilon}
\def\f{\varphi}

\def\s{\sigma}
\def\t{\theta}
\def\l{\lambda}
\def\LEQ{\leqslant}
\def\GEQ{\geqslant}

\def\Log{{\rm Log\,}}

\def\DST{\displaystyle}
\def\pr{{\bf Proof:\\ }}
\def\HF{\hfill{$\diamondsuit$\\ }}

\def\umLtbt{(1-\Log t)^{\beta_\theta}}
\def\fep{f_*^p(s)ds}
\def\dtt{\dfrac{dt}t}
\def\intH{H(x)dx}
\def\NN{\nonumber}

\def\up{{\frac1p}}
\def\ur{{\frac1r}}
\def\bt{{\beta_\theta}}
\def\supztu{{\displaystyle\sup_{0<t<1}}}

\newcommand{\umL}[3]{(1-{\rm Log\,}{#1})^{\frac{#2}{#3}}}
\newcommand{\suppp}[2]{\displaystyle\mathop{\rm sup\,}_{0<{#1}<{#2}}}
\newcommand{\K}[2]{K(f,t;L^{#1},L^{#2})}
\newcommand{\ummL}[5]{(1-{\rm Log\,}{#1})^{\frac{#2}{#3}+\frac{#4}{#5}}}
\def\Ztr{{Z_{\theta,r}}}

\title{\bf Characterization of interpolation between \\Grand, small or classical Lebesgue spaces}
\author{{\bf Alberto Fiorenza}\\
{\small Dipartimento di Architettura,
Universit\`a di Napoli "Federico II",}
{\small  via Monteoliveto, 3, I-80134 Napoli,  ITALY,}\\
{\small and Istituto per le Applicazioni del Calcolo "Mauro Picone" }\\
{\small Consiglio Nazionale delle Ricerche} 
{\small via Pietro Castellino, 111 I-80131Napoli,  ITALY}\\
{\small   e-mail: fiorenza@unina.it}\\ \ \\
{\bf Maria Rosaria Formica}\\
{\small Universit\`a degli Studi di Napoli "Parthenope",
via  Generale Parisi 13, 80132, Napoli, ITALY}\\
{\small   e-mail: mara.formica@uniparthenope.it}\\
 \ \\
 {\bf Amiran Gogatishvili}\\
 {\small  Institute of  Mathematics of the   Czech Academy of Sciences -
\v Zitn\'a, 115 67 Prague 1, CZECH REPUBLIC}\\ 
 {\small e-mail: gogatish@math.cas.cz}\\
\ \\
{\bf Tengiz Kopaliani }\\
{\small Faculty of Exact and Natural Sciences}\\
{\small Javakhishvili Tbilisi State University - University St. 2  0143 Tbilisi, GEORGIA}\\
{\small e-mail: tengiz.kopaliani@tsu.ge}\\
\ \\
{\bf Jean Michel Rakotoson\footnote{corresponding author: {\small rako@math.univ-poitiers.fr, jean.michel.rakotoson@univ-poitiers.fr}}}\\
{\small Laboratoire de Mathématiques et Applications - Université de Poitiers,}\\ {\small Avenue Marie et Pierre Curie,T\'el\'eport 2,}
{\small  BP 30179,86692 Futuroscope Chasseneuil Cedex, FRANCE} \\
{\small  rako@math.univ-poitiers.fr   jean.michel.rakotoson@univ-poitiers.fr }}
\date{\today}                                           

\def\G{\Gamma}
\def\O{\Omega}

\def\a{\alpha}
\def\b{\beta}

\def\f{\varphi}

\def\s{\sigma}
\def\t{\theta}
\def\LEQ{\leqslant}
\def\GEQ{\geqslant}

\def\Log{{\rm Log\,}}

\def\DST{\displaystyle}

\def\e{{e^{1-\frac1t}}}
\def\up{\frac1p}
\def\dss{\dfrac{ds}s}
\def\xd{\left(\dfrac x2\right)}

\def\R{{\rm I\!R}}
\begin{document}
\large
\maketitle
\begin{abstract}
In this paper, we show that the interpolation spaces between Grand, small or classical Lebesgue are so called Lorentz-Zygmund spaces or more generally  $G\G$-spaces. As a direct consequence of our results any Lorentz-Zygmund space
$L^{a,r}(\Log L)^\beta$, is an interpolation space in the sense of Peetre between either  two Grand Lebesgue spaces
or between two small spaces  provided that $ 1<a<\infty, \beta \not= 0$. The method consists in computing the so called K-functional of the interpolation space and in identifying the associated norm.

\end{abstract}

\tableofcontents
\section{Introduction. Main results}
Let $(X_0,||\cdot||_0),\ (X_1,||\cdot||_1)$ two Banach spaces contained continuously in a Hausdorff topological vector space (that is $(X_0,X_1)$ is a compatible couple).\\
For $g\in X_0+X_1,\ t>0$ one defines the so called $K$ functional $K(g,t;X_0,X_1)\dot=K(g,t)$ by setting
\begin{equation}\label{eq1}
K(g,t)=\inf_{g=g_0+g_1}\big(||g_0||_0+t||g_1||_1\big).
\end{equation}
For $0\LEQ\theta\LEQ1,\ 1\LEQ p\LEQ+\infty,\ \a\in\R$ we shall consider
$$(X_0,X_1)_{\theta,p;\a}=\Big\{g\in X_0+X_1,\ ||g||_{\theta,p;\a}=||t^{-\theta-\frac1p}\big(1-\Log t\big)^\a K(g,t)||_{L^p(0,1)}\hbox{ is finite}
\Big\}.$$
Here $||\cdot||_V$ denotes the norm in a Banach space $V$. The weighted Lebesgue space $L^p(0,1;\omega)$, $0< p\LEQ+\infty$ is endowed with the usual norm or quasi norm, where $\omega$ is a weight function on $(0,1)$.\\
Our definition of the  interpolation space is different from the usual one (see \cite{Bennett-Sharpley, Tartar}) since we restrict the norms on the interval $(0,1)$.\\ If  we consider ordered couple, i.e.  $X_1\hookrightarrow X_0$ and 
 $\a=0,\ (X_0,X_1)_{\theta,p;0}=(X_0,X_1)_{\theta,p}$ is the interpolation space as it is defined by  J. Peetre (see \cite{Bennett-Sharpley, Tartar, Bergh-Lofstrom}). \\


A particular attention will be brought to the so-called Grand Lebesgue space $L^{p),\a}(\O)$, with $\O$ a bounded (open) set of $\R^n$ whose measure is 1, $1<p<+\infty,\ \a>0$ defined as 
$$L^{p),\a}(\O)\!\!=\!\Big\{f:\O\to\R\hbox{\,measurable,\,}||f||_{p),\a}\!\!=\!\sup_{0<t<1}(1-\Log t)^{-\frac\a p}\left(\int_t^1f_*^p(\s)d\s\right)^{\frac1p}<+\infty
\Big\}$$
and its associated spaces $L^{(p',\a}(\O),\ \dfrac1{p'}+\dfrac1p=1$ defined as (see \cite{dff})

$$L^{(p',\a}(\O)\!\!=\left\{
f:\O\!\to\!\R\hbox{\,measurable,\,}
||f||_{(p',\a}\!\!=\!\!\int_0^1\!(1\!\!-\Log t)^{-\frac\a {p'}\!+\a\!-\!1}\!\!\left(\int_0^t\!\!f_*^{p'}(\s)d\s\right)^{\frac1{p'}}\!\dfrac{dt}t\!\!<\!\!+\!\infty
\right\}.$$
Here, $f_*$ is the decreasing rearrangement of $|f|$, say it is the generalized inverse of the distribution function
$$D_f(t)=measure \{x\in \O,  |f(x)|>t \}, t \in \R_+.$$
Many works related to Grand and small Lebesgue spaces have been recently done
(see for instance \cite{Ana, dff, FFR, fk, CobosKuhn}). 

We will show in particular, the

\begin{theo}\label{t4}\ \\
Let $1<p<q,\ \a>0$. Then
$$L^{q),\a}=(L^{p),\a},L^q)_{1,\infty;-\frac\a q}.$$
\end{theo}
An explicit equivalent of $K(f,t;L^{p),\a},L^q)$ is given in Theorem \ref{t3}.\\

For convenience, we will sometimes drop the dependence with respect to the domain $\O$ or $(0,1)$ and we shall write $L^p(\O)=L^p$\ldots etc\\ More, we will write sometimes $\DST\int_E f^m_*(t)dt=\int_Ef^m_*$.

\begin{defi}\label{d1}{(\bf Lorentz-Zygmund space)}\\
For $1\LEQ p,\ q\LEQ\infty,\ -\infty<\a<+\infty$, the Lorentz-Zygmund space $L^{p,q}(\Log\ L)^\a$ consists of all functions $f$ measurable such that 
$$||f||_{p,q;\a}=\begin{cases}\DST\left(\int_0^1\Big[t^{\frac 1p-\frac1q}(1-\Log t)^\a f_*(t)\Big]^qdt\right)^{\frac1q}&if\ 1\LEQ q<+\infty,\\\ \\
\DST\sup_{0<t<1}t^{\frac1p}(1-\Log t)^\a f_*(t)&if\ q=+\infty\end{cases}\hbox{ is finite.}$$
\end{defi}
One of the major theorems of the first section will be
\begin{theo}\label{t2}\ \\
Let $0<\theta<1,\ 1\LEQ r<+\infty$, $\a>0$, $1<p<q$. Then 
$$\Big(L^{p),\a},L^{q),\a}\Big)_{\theta,r}=L^{p_\theta,r}\big(\Log L\big)^{-\frac{\a}{p_\theta}}\qquad\hbox{
where }\DST\frac1{p_\theta}=\frac {1-\theta}p+\frac\theta q.$$
\end{theo}

We  will also use  the following extension of Generalized Gamma space (see \cite{FR2}).
\begin{defi}{\bf of a Generalized Gamma space with double weights}\\
Let $w_1,\ w_2$  be two weights on $(0,1)$, $m\in[1,+\infty]$, $1\LEQ p<+\infty$. We assume the following conditions:
\begin{itemize}
\item[c1)] There exists $K_{12}>0$ such that $w_2(2t)\LEQ K_{12}w_2(t)\ \forall\,t\in(0,1/2)$. The space
$L^p(0,1;w_2)$ is continuously embedded in $L^1(0,1)$.
\item[c2)] The function $\DST\int_0^tw_2(\s)d\s$ belongs to $L^{\frac mp}(0,1;w_1)$.
\end{itemize}
A generalized Gamma space with double weights is the set :
$$G\G(p,m;w_1,w_2)=\left\{v:\O\to\R\hbox{ measurable }\int_0^tv_*^p(\s)w_2(\s)d\s \hbox{ is in }L^{\frac mp}(0,1;w_1)\right\}.$$
\end{defi}

\begin{propri}\ \\
Let $G\G(p,m;w_1,w_2)$ be a Generalized Gamma space with double weights and let us define for $v\in G\G(p,m;w_1,w_2)$
$$\rho(v)=\left[\int_0^1w_1(t)\left(\int_0^tv_*^p(\s)w_2(\s)d\s\right)^{\frac mp}dt\right]^{\frac1m}$$
with the obvious change for $m=+\infty$.\\
Then,
\begin{enumerate}
\item $\rho$ is a quasinorm.
\item$G\G(p,m;w_1,w_2)$ endowed  with $\rho$ is a quasi-Banach function space.
\item If $w_2=1$ $$G\G(p,m;w_1,1)=G\Gamma(p,m;w_1).$$
\end{enumerate}
\end{propri}
\pr
\begin{enumerate}
\item 
Due to the property of the monotone rearrangement we have $\rho(v)=0$ if and only if $v=0$ and for $\lambda\in\R$ 
$\rho(\lambda v)=|\lambda|\rho(v)$. Let us show that
$$\rho(v_1+v_2)\LEQ (2K_{12})^{\frac1p}\Big(\rho(v_1)+\rho(v_2)\Big) \hbox{ for $v_1$ and $v_2$ in } G\G(p,m;w_1,w_2).$$
We have for all $\s\in(0,1),\ v=v_1+v_2$
$$(v_1+v_2)_*(\s)\LEQ v_{1*}\left(\dfrac\s2\right)+v_{2*}\left(\dfrac\s2\right).$$
Therefore, we have (using the triangle inequality)
\begin{eqnarray}
\left(\int_0^tv_*^p(\s)w_2(\s)d\s\right)^{\frac1p}
&\LEQ&\left(\int_0^tv_{1*}^p\left(\dfrac\s2\right)w_2(\s)d\s\right)^{\frac1p}+\left(\int_0^tv_{2*}^p\left(\dfrac\s2\right)w_2(\s)d\s\right)^{\frac1p}\label{eq15000}\\
&\LEQ&2^{\frac1p}\left(\int_0^{\frac t2}v_{1*}^p(s)w_2(2s)ds\right)^{\frac1p}+2^{\frac1p}\left(\int_0^{\frac t2}v_{2*}^p(s)w_2(2s)ds\right)^{\frac1p}\nonumber\\
&&\label{eq15001}
\end{eqnarray}
Using condition c1), we have
$$\left(\int_0^tv_*^p(\s)w_2(\s)d\s\right)^{\frac1p}
\LEQ(2K_{12})^{\frac1p}\left[\left(\int_0^tv_{1*}^p(s)w_2(s)ds\right)^{\frac1p}+\left(\int_0^tv_{2*}^p(s)w_2(s)ds\right)^{\frac1p}\right].$$
Using again the triangle inequality with the space $L^m(0,1;w_1)$ we derive that:
$$\rho(v_1+v_2)\LEQ(2K_{12})^{\frac1p}\Big[\rho(v_1)+\rho(v_2)\Big].$$

\item The function $\rho$ satisfies:\\
 If $0\LEQ v_1\LEQ v_2$ then $v_{1*}\LEQ v_{2*}$ everywhere on $(0,1)$ so that $$\rho(v_{1})\LEQ \rho(v_{2}).$$

If $0\LEQ v_{1k}\nearrow v$ almost-everywhere then by the Beppo-Levi's theorem we have for all $t\in[0,1]$
$$\lim_{k\to+\infty}\left(\int_0^tv^p_{1k*}(\s)w_2(\s)d\s\right)^{\frac1p}=\left(\int_0^t v_*^p(\s)w_2(\s)d\s\right)^{\frac1p}, $$
since $v_{1k*}\nearrow v_*$ everywhere on $[0,1]$.\\

Because $L^m(0,1;w_1)$ is a Banach function space, we deduce that
$$\lim_{k\to+\infty}\rho(v_{1k})=\rho(v).$$
The condition c2) implies that $$\rho(\chi_E)\LEQ\rho(1)<+\infty$$
for any $E\subset\O,\ \chi_E$ denoting its characteristic function.\\

To conclude that it is a complete space, we shall prove the  inequality in next Proposition \ref{piGG}, which has been already given in the frame of $G\G(p,m;w_1)$ \Big(with one weight (see \cite{FR1, FR2, FFR})\Big).\ \HF
\end{enumerate}

\begin{propo}{\bf Inequality for $G\G$}\label{piGG}\ \\
Let $v\in G\G(p,m;w_1,w_2),\ E\subset\O,\ |E|>0$. Then
$$\dfrac{\DST \rho(v)\left[\int_0^{|E|}w_2(\s)d\s\right]^{\frac1p}}{\DST\left[\int_0^{|E|}w_1(t)\left(\int_0^tw_2(s)ds\right)^{\frac mp}dt\right]^{\frac1m}}\GEQ
\left[\int_0^{|E|}v_*^p(\s)w_2(\s)d\s\right]^{\frac1p}. $$
\end{propo}
\pr
Let us set $$V(t)=\dfrac{\DST\int_0^t  v_*^p(\s) w_2(\s)d\s}{\DST\int_0^tw_2(\s)d\s}.$$
Then $V$ is decreasing since $v_*^p$ is decreasing. Therefore, we have for $1\LEQ m<+\infty$
$$\rho(v)=\left[\int_0^1 V^{\frac mp}(t)w_1(t)\left(\int_0^t w_2(\s)d\s\right)^{\frac mp}dt\right]^{\frac1m}
\GEQ V^{\frac1p}(|E|)\left[\int_0^{|E|}w_1(t)\left(\int_0^tw_2(\s)d\s\right)^{\frac mp}dt\right]^{\frac1m}.$$
This gives the inequality $m<+\infty$. \\
For $m=+\infty$, the argument is the same, since
\begin{eqnarray*}
\sup_{t\LEQ1}w_1(t)\left(\int_0^tv_*^p(\s)w_2(\s)d\s\right)^{\frac1p}&\GEQ&\sup_{t\LEQ|E|}V^{\frac1p}(t)w_1(t)\left(\int_0^tw_2(\s)d\s\right)^{\frac1p}\\
&\GEQ& V^{\frac1p}(|E|)\sup_{t\LEQ |E|}w_1(t)\left(\int_0^tw_2(\s)d\s\right)^{\frac1p}.
\end{eqnarray*}
\
\HF
One of our results which motivates the introduction of Generalized Gamma space with double weights is
\begin{theo}\label{tFund1}\ \\
Let $1<p<+\infty,\ 0<\t<1,\ 1\LEQ r<+\infty$. Then
$$(L^{p)},L^{(p})_{\t,r}=G\G(p,r;w_1,w_2)\hbox{ with } w_1(t)=t^{-1}(1-\Log t)^{\t r-1},\ w_2(t)=(1-\Log t)^{-1},\ t\in(0,1).$$
\end{theo}
The proof is given in the last theorem of this paper.
\section{Notations. Preliminary Lemmas}
For two positive quantities $A$ and $B$ depending on some parameters like functions \big(real number  as a function  on  $(0,1)$\big) we shall write $A\lesssim B$ if there exists a constant $c>0$  independent of the parameters such that $A\LEQ cB$, and $A\approx B$ if $B\lesssim A$ and $A\lesssim B$.

\begin{propo}\label{p1000}\ \\
Let $a>0,\ b\in\R,\ \Phi(t)=t^a(1-\Log t)^b$, $t\in(0,1)$. Then, there exists an invertible function $\f$ from $[0,1]$ into $[0,1]$ such that 
\begin{enumerate}
\item $t\approx\Phi(\f(t)) =\f(t)^a(1-\Log\f(t))^b$
\item $1+|\Log\f(t)|\approx 1+|\Log t|$ for $t\in(0,1)$.
\end{enumerate}
\end{propo}
\pr
Let us set $t_0=\begin{cases}e^{\frac {a-b}a}&if\ a<b,\\1&if\ a\GEQ b.\end{cases}$
Then $\Phi$ is strictly increasing on $[0,t_0]$. \\ 

Define $g(t)=\dfrac{\Phi(t \,t_0) }{\Phi(t_0)}$  for $t\in[0,1]$. Then $g$ is continuous and strictly increasing from $[0,1]$ into itself, and  $g(t)\approx t^a(1-\Log t)^b$ and then $1-\Log g(t)\approx 1-\Log t$. Setting $\f(t)=g^{-1}(t)$ the inverse of $g$  we have the result.\HF  
   \begin{propo}\label{p2}\ \\
 Let $\b\in\R,\ -\infty<\a<1$. Then, there exists $c_{\a\b}>0$:
 $$\int_{0}^at^{-\a}(1-\Log t)^{\b}dt\LEQ c_{\a\b}a^{1-\a}(1-\Log a)^{\b}\qquad\forall\,a\in[\,0,1\,].$$
 Moreover, if $\beta \GEQ 0$ then
$$ \int_{0}^at^{-\a}(1-\Log t)^{\b}dt\GEQ \frac1{1-\a}a^{1-\a}(1-\Log a)^{\b}.$$
 \end{propo}
\pr We start with the case $\beta>0$.\\
Let $k\in{\N}$ such $k-1\LEQ \b<k$. Then by integration by parts, we have
\begin{eqnarray*}
\int_0^{a}t^{-\a}(1-\Log t)^\b dt&=&c_{\a\b k}\int_0^at^{-\a}(1-\Log t)^{\b-k}dt+a^{1-\a}\sum_{j=0}^{k-1}c_{\a\b j}(1-\Log a)^{\b-j}\\
&&\hbox{\small (using the fact that  the function $t\rightarrow(1-\Log t)^{\b-k}$ is increasing).}\\
&\LEQ& c_{\a\b k}(1-\Log a)^{\b-k}\int_0^at^{-\a}dt+c_{\a\b k}a^{1-\a}(1-\Log a)^\b\\
&\LEQ& c_{\a\b}a^{1-\a}(1-\log a)^\b.\end{eqnarray*}

If $\b<0$ the inequality is still true since the function $ t\rightarrow(1-\Log t)^\b$ is increasing, we argue as in the last line
of the above proof. If $\b \GEQ 0$ then $ t\rightarrow(1-\Log t)^\b$ is decreasing, then the last result follows directly.
\ \HF
Using an argument of \cite{BR}, we have:
\begin{propo}\label{p2.3}\ \\
For any $\a<-1,\ \b\in\R$, we have
$$\int_a^1t^\a(1-\Log t)^\b dt\lesssim a^{\a+1}(1-\Log a)^\b.$$
\end{propo}
\pr
Let $\eps>0$ such that $\a+\eps<-1$. For  $N=1+\dfrac{|\b|}\eps,\ t^{-\eps }(N-\Log t)^\b$ is decreasing. So $\DST\int_a^1t^\a(1-\Log t)^\b dt\lesssim a^{-\eps}(N-\Log a)^\b\int_a^1t^{\a+\eps}dt\lesssim(1-\Log a)^\b a^{\a+1}$.\HF\\

For convenience,
 we recall some of those Hardy type inequalities
 (see  \cite{BR} Theorem 6.4, \cite{GOT} Lemma 2.7, Corollary~2.9)  that we will use.

\begin{theo}\label{t5.3}\label{BN}\ \\
Suppose  $\l>0,\ 1\LEQ b\LEQ +\infty$, and $-\infty<\b<+\infty$, $\Phi$ a nonnegative measurable function on $(0,1)$. Then
$$\int_0^1\left[t^{-\l}(1-\Log t)^\b\int_0^t\Phi(s)ds\right]^b\dfrac{dt}t\LEQ c\int_0^1\Big[t^{1-\l}(1-\Log t)^\b\Phi(t)\Big]^b\dfrac{dt}t,$$
and
$$\int_0^1\left[t^\l(1-\Log t)^\b\int_t^1\Phi(s)ds\right]^b\dfrac{dt}t\LEQ c\int_0^1\left[t^{1+\l}(1-\Log t)^\b\Phi(t)\right]^b\frac{dt}t.$$
The constant $c$ is independent of $\Phi$.\\
If $\Phi(t)=t^{\mu-1}\Phi_1(t),\ \mu>0,\ \Phi_1$ decreasing, then the above inequalities hold true for $0<b<1$.
\end{theo}
In the above formula when $b=+\infty$, the integral is replaced by the supremum. 
In the case where $\l=0$ we use the following Hardy inequalities given in Bennett-Rudnick (\cite{BR}~Theorem~6.5)
\begin{theo}\label{tBR6.5}\ \\
Suppose $1\LEQ a\LEQ \infty,\ \a\in\R$ and $\a+\dfrac1a \neq0$. Let $\psi $ be a nonnegative measurable function on $(0,1)$. Then, if $\a+\dfrac1a>0$,
\begin{equation}\label{eqBR6.13}
\left(\int_0^1\left[(1-\Log t)^\a\int_0^t\psi(s)ds\right]^a\dtt\right]^{\frac1a}
\LEQ c\left(\int_0^1\Big[t(1-\Log t)^{1+\a}\psi(t)\Big]^a\dtt\right]^{\frac1a}
\end{equation}
and $\a+\dfrac1a<0$,
\begin{equation}\label{eqBR6.14}
\left(\int_0^1\left[(1-\Log t)^\a\int_t^1\psi(s)ds\right]^a\dtt\right)^{\frac1a}
\LEQ c \left(\int_0^1\Big[t(1-\Log t)^{1+\a}\psi(t)\Big]^a\dtt\right)^{\frac1a}.
\end{equation}
The constant $c$ is independent of $\psi$.
\end{theo}

Let $(\RR,\mm) $ be a measure space and $\M(\RR,\mm)$ be the set of all $\mu$ measurable functions over $\RR$ .
A~Banach space $X=X(\RR,\mm)$ of $\mu$-measurable complex-valued functions in $\M(\RR,\mm)( \hbox { set of all $\mu$ measurable functions 
over $\RR$ } )$, equipped with the
norm $\|\cdot\|_{X}$, is said to be a {\it
rearrangement-invariant Banach function space} (shortly \textit{r.i.~space}) over $(\RR,\mu)$ (or over $\RR$ with respect to $\mu$) if the following five axioms hold:
\begin{itemize}
\item[(P1)]\qquad $0\le g \le f$ $\mu$-a.e.\  implies
$\|g\|_{X} \le \|f\|_{X}$; 
\item[(P2)]\qquad $0\le f_n
\nearrow f$ $\mu$-a.e.\ implies $\|f_n\|_{X} \nearrow
\|f\|_{X}$; 
\item[(P3)]\qquad $\|\chi_E\|_{X}<\infty$ for every $E\subset\RR$ of finite measure;
\item[(P4)]\qquad for every $E\subset\RR$ with $\mu(E)<\infty$ there exists a~constant $C_E$ such that $\int_{E}
|f(x)|\,d\mu(x) \le C_E \|f\|_{X}$ for every $f\in X$;
\item[(P5)]\qquad $\|f\|_{X} = \|g\|_{X}$ whenever $f_* = g_ *$.
\end{itemize}

Given an~r.i.~space $X$ on $(\RR,\mu)$, the set
$$
X'=\left\{f\in\M(\RR,\mu):\,\int_{\RR}|f(x)g(x)|\,d\mu(x)<\infty\
\textup{for every } g\in X\right\},
$$
equipped with the norm
$$
\|f\|_{X'}=\sup_{\|g\|_{X}\leq1} \int_{\RR}|f(x)g(x)|\,d\mu(x)
=\sup_{\|g\|_{X}\leq1} \int_{0}^{\mu(\RR)}f_*(t)g_*(t)\,dt,
$$
is called the \textit{associate space} of $X$. It turns out
that $X'$ is again an~r.i.~space over $\RR$ with respect to $\mu$ and that $X''=X$.

For every r.i.~space $X$ over $(\RR,\mu)$, there exists a~unique
r.i.~space $\overline X$ over $(0,\infty)$ with respect to the
one-dimensional Lebesgue measure, satisfying
$$
\|f\|_{X}=\|f_*\|_{\overline X}
$$
for every $f\in X$. This space, equipped with the norm
$$
\|f\|_{\overline X}=\sup_{\|g\|_{X'}\leq1}
\int_0^{\mu(\RR)}f_*(t)g_*(t)\,dt,
$$
is called the \textit{representation space} of $X$. 

The fundamental function of a r.i. Banach space,  $X$ , is defined
by 
\[\Phi_X(t) = \|\chi_{[0,t]}\|_X, \quad  t\in (0,\mu(\RR)).
\].
There is no loss of generality if we assume $\Phi_X $ to be
positive, nondecreasing, absolutely continuous far from the origin, concave and to
satisfy 
\[\Phi_X(t)\Phi_{X'} (t) = t  \quad \text{for all} \quad    t\in (0, \mu(\RR)).\]
Details and further material on r.i.~spaces can be found in~\cite[Chapter~2]{Bennett-Sharpley}.

Note that  (see \cite{CF}) grand Lebesgue spaces $L^{p),\a}(\O)$ and small Lebesgue spaces $L^{(p,\a}(\O)$ are r.i. spaces over $\O$ and 
 \[\Phi_{L^{p),\a}}(t)\approx t^{\frac1p}\umL t {-\a}p\quad \text{and}\quad  \Phi_{L^{(p,\a}}(t)\approx t^{\frac1p}\umL t {\a}{p'}
\]
 $(L^{p),\a})'= L^{(p',\a}$. 

Let $A_0$  and $A_1$ be two r.i. Banach spaces  over $[0,1]$ and let $\Phi_0,\, \Phi_1$  be respectively their fundamental
functions.  We suppose that the following conditions are satisfied:
\begin{description}
	\item{i)} There exists a constant C such that, for $i = 0; 1$ and for all $t > 0$
\begin{equation} \int_0^t \frac{ds}{\Phi_i(s)}\le \frac{Ct}{\Phi_i(t)},  \tag{C.0} \label{C0}
\end{equation}  
\item{ii)} There exists a constant C such that, for all $t >0$  
\begin{equation}\frac{\Phi_1(t)}{\Phi_0(t)}\left \|\frac{\chi_{[0,t]}}{\Phi_1}\right\|_{A_0}\le C, \tag{C.1} \label{C1}
\end{equation}
\begin{equation}
 \frac{\Phi_0(t)}{\Phi_1(t)} \left\|\frac{\chi_{[t,1]}}{\Phi_0}\right\|_{A_1}\le C.  \tag{C.2} \label{C2}
\end{equation}
\end{description}

\begin{lem}{\cite[Lemma 2.2]{BRF}}  \label{lemBR} Let $A_0$ and $A_1$ be two r.i. Banach function spaces satisfied   \eqref{C0}, \eqref{C1} and \eqref{C2} 
Then 
\[ K\left( f,\frac{\Phi_0(t)}{\Phi_1(t)}; A_0,A_ 1\right)\approx  \left \|f_*\chi_{[0,t]}\right\|_{A_0} + \frac{\Phi_0(t)}{\Phi_1(t)}\left \|f_*\chi_{[t,1]}\right\|_{A_1},\quad  0<t<1.
\] 
\end{lem}

{\bf Remarks on the choice of the method}\\
We have chosen a direct method  for the proof of our results by computing the $K$-functionals. In some part of the manuscript (for instance Theorem 5.1) we can adopt an alternative proof as using limiting  reiteration theorems \cite{AEEK, CC, Ew, FS}. Although, we observe that it is not possible to get  the result  without computations. The feature of our method is that, as a byproduct, we make explicit the behavior of the $K$-functional.

\section{Computation of some $K$-functionals and characterization of the interpolation spaces
 $(L^{p),\a},L^{q),\a})_{\theta,r}$}

  We shall need few lemmas before reaching the proof of Theorem \ref{t4} and Theorem \ref{t2}.

\begin{lem}\label{l1}\ \\
Let $\s=\dfrac p{p-1},\ 1<p<+\infty$, $\a>0$. Then
$$\sup_{0<t<1}t^{-1}(1-\Log t)^{-\frac\a p}\int_0^{t^\s}f_*(s)ds\LEQ c_{\a p}\sup_{0<t<1}(1-\Log t)^{-\frac\a p}
\left(\int_t^1f^p_*(s)ds\right)^{\frac1p}$$
$0<t<1$, for some constant $c_{\a p}>0$.
\end{lem}
\pr
We have
\begin{eqnarray}
 N_1&\equiv&\sup_{0<t<1}t^{-1}(1-\Log t)^{-\frac\a p}\int_0^{t^{\s}}f_*(s)ds\nonumber\\
 &=&\sup_{0<t<1}t^{-1}(1-\Log t)^{-\frac\a p}\int_0^{t^\s}\left(\dfrac2s\right)^{\frac1p}\left(\int_{\frac s2}^sd\tau\right)^{\frac1p}f_*(s)ds\nonumber\\
 &\lesssim&\sup_{0<t<1}t^{-1}(1-\Log t)^{-\frac\a p}\int_0^{t^\s}s^{-\frac1p}\left(\int_{s/2}^sf^p_*(\tau)d\tau\right)^{\frac1p}ds\nonumber\\
 &\lesssim&\sup_{0<t<1}t^{-1}\left(1-\Log t\right)^{-\frac\a p}\left(\int_0^{t^\s}s^{-\frac1p}\left(1-\Log\dfrac s2\right)^{\frac\a p}ds\right)\NN\\&&\times
 \sup_{0<s<t}\left(1-\Log\dfrac s2\right)^{-\frac\a p}\left(\int_{s/2}^1 f^p_*(\tau)d\tau\right)^{\frac1p}\nonumber\\
 &&\ \label{eqx1}
\end{eqnarray}
Applying Proposition \ref{p2}, we have
\begin{equation}\label{eqx2}
\int_0^{t^\s}s^{-\frac1p}\left(1-\Log\dfrac s2\right)^{\frac\a p}ds\LEQ c_{\a p} t^{\s(1-\frac1p)}(1-\Log t)^{\frac \a p}
\end{equation}
From the two last inequalities, we have:
$$N_1\LEQ c_{\a p}\sup_{0<t<1}(1-\Log t)^{-\frac \a p}\left(\int_t^1f_*^p\right)^{\frac1p}\approx ||f||_{p),\a}.$$
\ \HF
The next lemma has been already proved  (see Theorem 4.1 \cite{fk}) for a characterization of grand Lebesgue spaces as interpolation spaces between Lebesgue spaces, obtained using the Holmstedt's formula as well. 
\begin{lem}\label{l2}\ \\
For $1<p<+\infty,\ \a>0$ we have
$$L^{p),\a}=(L^1,L^p)_{1,\infty;-\frac\a p}.$$
\end{lem}
\pr
According to the Holmstedt's formula (see \cite{Bennett-Sharpley}), we have for all $f\in L^1+L^p$, and for all $t\in(0,1)$
$$K(f,t;L^1,L^p)\approx\int_0^{t^\s}f_*(s)ds+t\left(\int_{t^{\s}}^1f^p_*(s)ds\right)^{\frac1p}\dot=K_{1p}(f,t)\quad\hbox{with }\s=\dfrac p{p-1}.$$
Making use of the definition of the norm in $(L^1,L^p)_{1,\infty;-\frac\a p}$
$$||f||_{1,\infty;-\frac\a p}=\sup_{0<t<1}t^{-1}(1-\Log t)^{-\frac\a p}K_{1p}(f,t)$$
$$\approx\sup_{0<t<1}t^{-1}(1-\Log t)^{-\frac\a p}\int_0^{t^\s}f_*(s)ds+\sup_{0<t<1}(1-\Log t)^{-\frac\a p}\left(\int_{t^\s}^1f_*^p(s)ds\right)^{\frac1p}.$$
From Lemma \ref{l1}, we deduce that
$$||f||_{1,\infty;-\frac\a p}\approx\sup_{0<t<1}(1-\Log t)^{-\frac\a p}\left(\int_t^1 f_*^p(s)ds\right)^{\frac1p}=||f||_{p),\a} .$$
Noticing $1-\Log t\approx 1-\Log t^{\frac 1\s}.$
\HF
\begin{lem}\label{l4}\ \\
Let $1<p<q,\ \a>0$. Then
\begin{equation}\label{eqx3}
\sup_{0<t<1}t^{\frac1q-\frac1p}(1-\Log t)^{-\frac\a q}\left(\int_0^t f_*^p(s)ds\right)^{\frac1p}
\LEQ c_{\a pq}\sup_{0<t<1}(1-\log t)^{-\frac\a q}\left(\int_t^1f_*^q(s)ds\right)^{\frac1q}
\end{equation}
for some constant $c_{\a pq}>0$.
\end{lem}
\pr
The relation (\ref{eqx3}) is equivalent to
$$\sup_{0<t<1}t^{-1}(1-\Log t)^{-\frac{\a p}q}\left(\int_0^{t^\s}f_*(s)ds\right)
\lesssim\sup_{0<t<1}(1-\Log t)^{-\frac{\a p}q}\left(\int_t^1f_*(s)^{\frac qp}ds\right)^{\frac pq}$$$\hbox{ with } \s=\dfrac r{r-1},\ r=\dfrac qp>1.$\\

So we deduce Lemma \ref{l4} from Lemma \ref{l1}, replacing $p$ by $r$.
\ \HF
\begin{theo}\label{t2bb}\ \\
Let $1\LEQ p<q,\ \a>0$. Then
$$L^{q),\a}(\O)=\Big(L^p(\O),L^q(\O)\big)_{1,\infty;-\frac\a q}\ .$$
\end{theo} 

{\bf Proof of Theorem \ref{t2bb}}\\
According to Holmstedt's formula \cite{Bennett-Sharpley} we have for all $f\in L^p+L^q$, all $t\in(0,1)$
\begin{equation}\label{eqx4}
K(f,t;L^p,L^q)\approx\left(\int_0^{t^\s}f_*^p(s)ds\right)^{\frac1p}+t\left(\int_{t^\s}^1f_*^q(s)ds\right)^{\frac1q}\dot=K_{pq}(f,t)
\end{equation}
with $\s=\dfrac{pq}{q-p}$.
\\According to the norm of $f$ in $(L^p,L^q)_{1,\infty;-\frac\a q}$ we have:
\begin{eqnarray}
||f||_{1,\infty;-\frac\a q}&\approx&\sup_{0<t<1}t^{-1}(1-\Log t)^{-\frac\a q}K_{pq}(f,t)\nonumber\\
&\approx&\sup_{0<t<1}t^{-1}(1-\Log t)^{-\frac\a q}\left(\int_0^{t^\s}f_*^p(s)ds\right)^{\frac1p}
+\sup_{0<t<1}(1-\Log t)^{-\frac\a q}\left(\int_{t^\s}^1f_*^q(s)ds\right)^{\frac1q}\NN\\\label{eqx5}\\
&=& I_1+I_2.\nonumber
\end{eqnarray}
The first term $I_1$ of relation (\ref{eqx5}) is equivalent to
$$\sup_{0<\tau<1}\tau^{\frac1q-\frac1p}(1-\log\tau)^{-\frac\a q}\left(\int_0^\tau f_*^p(s)ds\right)^{\frac1p}.$$
(making use of the change of variables $\tau=t^\s$ and knowing as for $I_2, $ $1-\Log\tau^{\frac1\s}\approx1-\Log\tau$).\\
According Lemma \ref{l4}, this last term is dominated by $I_2$.\\
Therefore we have
$$||f||_{1,\infty;-\frac\a q}\approx I_1+I_2\approx I_2=||f||_{q),\a}.$$
\ \HF
 
\begin{theo}\label{t3}{\bf (Computation of the $K$-functional for the couple $(L^{p),\a},L^q)$)}\ \\
Let $1<p<q,\ \a>0$. Then
\begin{eqnarray*}K(f,t;L^{p),\a},L^q)&\approx&
\sup_{0<s<\f(t)}
(1-\Log s)^{-\frac\a p}
\left(\int_s^{\f(t)}f_*^p(x)dx\right)^{\frac1p}
+t\left(\int_{\f(t)}^1f_*^q(s)ds\right)^{\frac1q}\\&=&\OV K_{pq}(f,t)\end{eqnarray*}
where $\f$ is the inverse of the increasing function $\psi(t)=t^{\frac1p-\frac1q}(1-\Log t)^{-\frac\a p},\quad t\in(0,1)$.
Thus $$t=\f(t)^{\frac1p-\frac1q}(1-\Log \f(t))^{-\frac{\a} p}.$$
\end{theo}
{\bf Proof of Theorem \ref{t3}}\\  To apply Lemma  ~\ref{lemBR}, we need to check  that the conditions \eqref{C0},  \eqref{C1} and \eqref{C2}  are satisfied. 
Let $\Phi_0(t)=\left\|\chi_{[0,t]}\right \|_{L^{p),\alpha}}$ and $\Phi_1(t)=\left\|\chi_{[0,t]}\right \|_{L^{q}}$.
As we have  
\[ \Phi_0(t)\approx t^{\frac{1}{p}}(1-\Log t)^{-\frac{\alpha}{p}}, \quad  \Phi_1(t)\approx t^{\frac{ 1}{q}},\] 
we have that $\psi(t)=\frac{\Phi_0(t)}{\Phi_1(t)}$..
The conditions \eqref{C0}  \eqref{C1} and \eqref{C2}  easily follow by using  the Proposition~\ref{p2} and Proposition~\ref{p2.3}. and using   Lemma~\ref{lemBR}
we get  desired  result.
  
\ \vspace{-2cm} \HF\\

By same way we can obtain following theorem 
\begin{theo}\label{t3.3}{\bf (Computation of the $K$-functional for the couple $(L^{p),\a},L^{q),\alpha})$)}\ \\
Let $1<p<q,\ \a>0$. Then
\begin{eqnarray*}K(f,t;L^{p),\a},L^{q),\alpha})&\approx&
\sup_{0<s<\f(t)}
(1-\Log s)^{-\frac\a p}
\left(\int_s^{\f(t)}f_*^p(x)dx\right)^{\frac1p}\\
&+&t \sup_{\f(t)<s<1}(1-\Log s)^{-\frac\a q}\left(\int_s^1f_*^q(x)dx\right)^{\frac1q},\end{eqnarray*}
where $\f$ is the inverse of the increasing function $\psi(t)=t^{\frac1p-\frac1q}(1-\Log t)^{-\frac\a p+\frac\a q},\quad t\in(0,1)$.
Thus $$t=\f(t)^{\frac1p-\frac1q}(1-\Log \f(t))^{-\frac{\a} p + \frac\a q}.$$
\end{theo}
{\bf Proof of Theorem \ref{t4}}\\
Let $f$ be in $(L^{p),\a},L^q)_{1,\infty;-\frac\a q}$. Then, its norm is 
$$||f||_{1,\infty;-\frac\a q}=\supztu t^{-1}\umL{t}{-\a}{q}\K{p),\a}q.$$
Following Theorem \ref{t3}, this expression gives:
\begin{eqnarray}\label{eqx30}
||f||_{1,\infty;-\frac\a q}&\approx&\suppp t{ 1} t^{-1}\umL t{-\a}q
\suppp s{\f(t)}\umL s{-\a}p
\left(\int_s^{\f(t)}f_*^p(\tau)d\tau\right)^{\frac1p}\nonumber\\&&+\suppp{ t}{1}\umL t{-\a}{q}\left(\int_{\f(t)}^1f_*^q(\tau)d\tau\right)^{\frac1q}=J_1+J_2.
\end{eqnarray}
But we have by the definition of $\f$,
$$t^{-1}=\f(t)^{\frac1q-\frac1p}(1+|\Log\f(t)|)^{\frac\a p},\quad1+|\Log\f|\approx1+|\Log t|.$$
Thus the first term $J_1$ of (\ref{eqx30}), after a change of variables, gives
\begin{equation}\label{eqx31}
J_1\approx\suppp t1 t^{\frac{p-q}{qp}}\umL t{\a(q-p)}{pq}\suppp s t\umL s{-\a}p
\left(\int_s^tf_*^p(\tau)d\tau\right)^{\frac1p}.
\end{equation}
But we can bound the last term of (\ref{eqx31}) as:
$$\suppp st\umL s{-\a}p\left(\int_s^tf_*^p(\tau)d\tau\right)^{\frac1p} 
\LEQ t^{\frac{q-p}{pq}}\suppp s t\umL s{\a(p-q)}{qp}\suppp st\umL s{-\a}q\left(\int_s^tf_*^q(\tau)d\tau\right)^{\frac1q}$$
(using Hölder inequality and introducing new factor).\\
Therefore, we can estimate $J_1$ (after simplifying):
\begin{equation}\label{eqx32}
J_1\lesssim\suppp s 1\umL s{-\a}q\left(\int_s^1f_*^q(\tau)d\tau\right)^{\frac1q}=||f||_{q),\a}.
\end{equation}
While for the second term $J_2$, we have:
\begin{equation}\label{eqx33}
J_2\approx\suppp t 1\umL{\f(t)}{-\a}q\left(\int_{\f(t)}^1f_*^q(\tau)d\tau\right)^{\frac1q}
\end{equation}
which implies
\begin{equation}\label{eqx34}
J_2\approx\suppp s1\umL s{-\a}q\left(\int_s^1f_*^q(\tau)d\tau\right)^{\frac1q}=||f||_{q),\a}.
\end{equation}
We then have from relations (\ref{eqx30}) to (\ref{eqx34})
\begin{equation}\label{eqx35}
J_2\lesssim
||f||_{1,\infty;-\frac\a q}\approx 
 J_1+J_2
\end{equation}
This shows the result.\HF
Next, we want to prove Theorem \ref{t2}.\\
For this, we shall need the:
\begin{lem}\label{lxxx}\ \\
Let $1<p<q,\ \a>0, \ 0<\theta<1$, $1\LEQ r<+\infty$. Then 
$$\Big(L^{p),\a},L^{q),\a}\Big)_{\theta,r}=\Big(L^{p),\a},L^q\Big)_{\theta,r;-\frac{\a \theta}q}.$$
\end{lem}
\pr
From Theorem \ref{t4}, we know that
$$\Big(L^{p),\a},L^{q),\a}\Big)_{\theta,r}=\Big(L^{p),\a},\big(L^{p),\a},L^q\big)_{1,\infty;-\frac{\a} q}\Big)_{\theta,r}.$$
In this case, we can apply the Holmstedt formula (see Theorem 3.1 in \cite{GOT}) valid  also for the extreme case $\theta_1=1$.
\begin{equation}\label{eqx36}
K\Big(f,t;L^{p),\a},(L^{p),\a},L^q)_{1,\infty;-\frac\a q}\Big)\dot \equiv K_{\a}(t)
\approx t\sup_{\psi(t)<s<1}s^{-1}\umL s{-\a}q \ K(f,s;L^{p),\a},L^q)
\end{equation}
where $\psi$ is an invertible  function from $[0,1]$ into itself according to Proposition \ref{p1000} 
so that 
$$t\approx \psi(t)\umL{\psi(t)}\a q\qquad\hbox{ for } t\in[\,0,1\,].$$

Considering the norm $||f||_{\theta,r}$ of $f$ in $\big(L^{p),\a},L^{q),\a}\big)_{\theta,r}$
we have:
\begin{eqnarray*}
||f||_{\theta,r}&\approx&\Big\|t^{-\theta}K_\a(t)\Big\|_{L^r(0,1;\frac{dt}t)}\\
&\approx&\Big\|t^{1-\theta}\sup_{\psi(t)<s<1}s^{-1}\umL s{-\a}q K(f,s;L^{p),\a},L^q)\Big\|_{L^r(0,1;\frac{dt}t)}\\
&\approx&\Big\|\psi(t)^{1-\theta}\umL {\psi(t)}{\a(1-\theta)}q\sup_{\psi(t)<s<1}s^{-1}\umL s{-\a}qK(f,s;L^{p),\a},L^q)\Big\|_{L^r(\frac{dt}t)}\\
&\approx&\Big\|t^{1-\theta}\umL t{\a(1-\theta)}q\sup_{t<s<1}s^{-1}\umL s {-\a}q K(f,s;L^{p),\a},L^q)\Big\|_{L^r(\frac{dt}t)}.
\end{eqnarray*}
The last line is obtained after the change of variables $t \rightarrow \psi(t)$. To control this last term, we need
\begin{propo}\label{pxxx}\ \\
Let $K_d$ be a decreasing nonnegative function on $(0,1)$ and $p,\ q,\ \a,\ \t,r$ are as in Lemma \ref{lxxx}. Then
$$I_r\dot\equiv\int_0^1\left[t^{1-\theta}\umL t{\a(1-\theta)}q\sup_{t<s<1}\umL s{-\a}qK_d(s)\right]^r\dfrac{dt}t$$
$$\approx\int_0^1\Big[t^{1-\theta}\umL t{-\a\theta}qK_d(t)\Big]^r\frac{dt}t\dot\equiv I_d.$$
 \end{propo}
 \pr
Since $\DST\sup_{t<s<1}\umL s{-\a}q K_d(s)\GEQ\umL t{-\a}qK_d(t),$
 we deduce
$$I_r\GEQ\int_0^1\Big[t^{1-\theta}\umL t{-\a \theta}qK_d(t)\Big]^r\dfrac{dt}t=I_d.$$
For the reverse inequality, since $K_d$ is decreasing and $s\to\umL s{-\a r}q$ increasing then,
\begin{equation}\label{eqx37}
\umL s {-\a r}qK^r_d(s)\lesssim\int_{s/2}^s\umL x{-\a r} qK_d^r(x)\dfrac{dx}x
\lesssim\umL{{\dfrac s2}}{-\a r}q K_d^r\left(\dfrac s2\right).
\end{equation}
Setting 
$$I_{1r}=\int_0^1t^{(1-\theta)r}\umL t {\a r(1-\theta)}q\sup_{t<s<1}
\left(\int_{s/2}^s\umL x{-\a r}q K_d^r(x)\dfrac{dx}x\right)\dfrac{dt}t,$$ 
thus, we can bound, the integral $I_r$ by $I_{1r}$:
\begin{equation}\label{eqx38}
I_r=\int_0^1t^{(1-\theta)r}\umL t{\a r(1-\theta)}q\sup_{t<s<1}\umL s{-\a r}qK_d^r\left({s}\right)\dfrac{dt}t
\lesssim I_{1r}
\end{equation}
and 
$$I_{1r}\lesssim\int_0^1t^{(1-\theta)r}\umL t{\a r(1-\theta}q\sup_{t<s<1}\umL{\dfrac s2}{-\a r}q K_d^r\left(\dfrac s2\right)\dfrac{dt}t=I_{2r}.$$
By a change of variables in the term $I_{2r}$, we have
\begin{equation}\label{eqx39}
I_{2r}\lesssim I_r.
\end{equation}
Therefore, relations (\ref{eqx38}) and (\ref{eqx39}) imply:
$$I_r\approx I_{1r}.$$

We may bound $I_{1r}$ as:
$$I_{1r}\LEQ\int_0^1t^{(1-\theta)r}\umL t{\a r(1-\theta)}q\left(\int_{t/2}^1\umL x{-\a r}qK^r_d(x)\dfrac{dx}x\right)\dfrac{dt}t.$$
By Fubini's Theorem, the upperbound in this inequality gives:
$$I_{1r}\lesssim\int_0^1\umL x{-\a r}qK_d^r(x)\left(\int_0^{2x}t^{(1-\theta)r}\umL t{\a r(1-\theta)} q\dfrac{dt}t\right)\dfrac{dx}x.$$
From Proposition \ref{p2}, we have
\begin{equation}\label{eqx40}
\int_0^{2x}t^{(1-\theta)r-1}\umL t{\a r(1-\theta)}q dt\LEQ c_{0q\a}x^{(1-\theta )r}\umL x{\a r(1-\theta)} q.
\end{equation}
From relation (\ref{eqx40}), the last estimate $I_{1r}$ becomes
$$I_{1r}\lesssim\int_0^1x^{(1-\theta)r}\umL x{-\a r\theta}qK^r_d(x)\dfrac{dx}x=I_d.$$
With relation (\ref{eqx38}) we get:
$$I_r\lesssim I_{1r}\lesssim I_d\LEQ I_r.$$
\ \HF

{\bf End of proof of Lemma \ref{lxxx}}\ \\
We apply Proposition  \ref{pxxx} with 
$$K_d(t)=t^{-1}\K{p),\a}q,\ t\in(0,1)$$ 
to derive the:
\begin{eqnarray*}
||f||_{\theta, r}&\approx&\Big\|t^{1-\theta}\umL t{\a(1-\theta)}q\sup_{t<s<1}\umL s{-\a}qK_d(s)\Big\|_{L^r(\frac{dt}t)}\\
&\approx&I^{\frac1r}_d=||f||_V\hbox{ with }V=\big(L^{p),\a},L^q\big)_{\theta,r;-\frac{\a\theta}q}.
\end{eqnarray*}
\ \HF

{\bf Proof of Theorem \ref{t2}}\\
Let 
$$V=\Big(L^{p),\a},L^{q),\a}\Big)_{\theta,r} \quad \text{and }\quad L^{p_\theta,r}(\Log L)^{-\frac{\a}{p_\theta}}=V_\theta$$ 
with
$\dfrac1{p_\theta}=\dfrac{1-\theta}p+\dfrac\theta q$.\\

We recall that, for $f\in V$
\begin{equation}\label{eqx50}
|f||_V^r=\int_0^1\left[t^{-\theta}\umL t{-\a\theta}q\K{p),\a}{q),\a}\right]^r\dfrac{dt}t.
\end{equation} 

Applying Theorem \ref{t3.3}, we have
\begin{eqnarray*}
||f||^r_V&\approx&\int_0^1\left[t^{-\theta}\sup_{0<s<\f(t)}\umL s{-\a}p\left(\int_s^{\f(t)}f_*^p(\tau)d\tau\right)^{\frac1p}\right]^r\dfrac{dt}t\\
&&+\int_0^1\left[t^{1-\theta}\sup_{\f(t)<s<1}\umL s{-\a}q\left(\int_{s}^1f_*^q(\tau)d\tau\right)^{\frac1q}\right]^r\dfrac{dt}t
\\&=&II_1+II_2.
\end{eqnarray*}

From the first term $II_1$, we know that
$$t=\f(t)^{\frac1p-\frac1q}\ummL{\f(t)}{-\a}p{\a}q$$
and $$1-\Log t\approx1-\Log\f(t).$$
Therefore, we have
\begin{equation}\label{eqx51}
t^{-\theta}\approx\f(t)^{\theta(\frac1q-\frac1p)}\umL{\f(t)}{-\a\theta(p-q)}{qp}.
\end{equation}
Replacing the integrand in the integral by the $II_1$ of (\ref{eqx51}) and making the change of variables $x=\f(t)$ or equivalently $t=\psi(x)$, we have
$$II_1\approx\int_0^1\left[x^{\theta(\frac1q-\frac1p)}\umL x{-\a\theta(p-q)}{qp}\sup_{0<s<x}
\umL s{-\a}p\left(\int_s^xf^p_*(\tau)d\tau\right)^{\frac1p}\right]^r\dfrac{dx}x.$$
Thus the monotony of $f_*$ leads to:
$$II_1\gtrsim\int_0^1\left[x^{\theta(\frac1q-\frac1p)}\umL x{-\a\theta(p-q)}{qp}\umL
{\dfrac x2}{-\a}px^{\frac1p}f_*(x)\right]^r\dfrac{dx}x$$
\begin{equation}\label{eq8000}
II_1\gtrsim\int_0^1\left[t^{\frac{1-\theta}p+\frac\theta q}\umL t{-\a(q(1-\theta)+\theta p)}{pq}
f_*(t)\right]^r\dfrac{dt}t=||f||^r_{V_\theta}.
\end{equation}
In particular, we have shown that
\begin{equation}\label{eqx56}
||f||_V\gtrsim||f||_{V_\theta}.
\end{equation}
For the  reverse of relation (\ref{eqx56}), let us consider
 $\eps>0$ small enough so that
$$\theta\left(\dfrac1q-\dfrac1p\right)+\dfrac\eps p<0.$$
$$II_1\lesssim\int_0^1\left[x^{\theta(\frac1q-\frac1p)}\umL x{-\a\theta(p-q)}{qp}\sup_{0<s<x}\umL s{-\a }p\left(\int_s^xt^{\eps-1}t^{1-\eps}f_*^p(t)dt\right)^{\frac1p}\right]^r\dfrac{dx}x.$$
Then
$$II_1\lesssim\int_0^1\left[x^{\theta(\frac1q-\frac1p)+\frac\eps p}\umL x{-\a}{p_\theta}\sup_{0<t<x}t^{\frac{1-\eps}p}f_*(t)\right]^r\dfrac{dx}x.$$
We have the following
\begin{propo}\label{pxx2}\ \\
For all $x>0$
$$\sup_{0<t<x}t^{\frac{1-\eps}p}f_*(t)\LEQ 2(\Log 2)^{\frac1{r'}}\left(\int_0^xs^{r\frac{1-\eps}p-1}f_*^r(s)ds\right)^{\frac1r},$$
with $\dfrac 1{r'}+\dfrac1r=1$.
\end{propo}
\pr
For any $s\in\left[\dfrac t2,t\right]$, we have
$$t^{\frac{1-\eps}p-1}f_*(t)\LEQ s^{\frac{1-\eps}p-1}f_*(s),$$
since $\dfrac1p-1-\dfrac\eps p<0$.\\
Integrating between $\dfrac t2$ and $t$, and using H\"older inequality, we have
$$t^{\frac{1-\eps}p}f_*(t)\LEQ2\left(\int_{t/2}^ts^{\frac{1-\eps}pr-1}f_*^r(s)ds\right)^{\frac1r}(\Log 2)^{\frac1{r'}}$$
$$t^{\frac{1-\eps}p}f_*(t)\LEQ 2\left(\int_0^t s^{\frac{1-\eps}pr-1}f^r_*(s)ds\right)^{\frac1r}(\Log2)^{\frac1{r'}}.$$
From which, we derive the result.\HF

Using this Proposition \ref{pxx2}, setting $m=\theta\left(\dfrac1q-\dfrac1p\right)+\dfrac\eps p <0$,
\begin{equation}\label{eqx52}II_1\lesssim\int_0^1x^{mr-1}\umL x{-\a r}{p_\theta}\left(\int_0^xs^{\frac{1-\eps}pr-1}f^r_*(s)ds\right)dx.
\end{equation}

Applying Hardy type inequality (see Theorem \ref{BN} and \cite{BR}, Theorem 6.4 ),  relation (\ref{eqx52}) becomes
\begin{equation}\label{eqx55}
II_1\lesssim\int_0^1\left[s^{\frac1{p_\theta}}\umL s{-\a}{p_\theta}f_*(s)\right]^r\dfrac{ds}s=||f||^r_{L^{p_\theta,r}(\Log L)^{-\frac{\a}{p_\theta}}}.
\end{equation}

It remains to show that the second term satisfies:
$$II_2=\int_0^1\left[t^{1-\theta}\sup{\f(t)<s<1}\umL s{-\a\theta} q\left(\int_{s}^1f_*^q(\tau)d\tau\right)^{\frac1q}\right]^r\dfrac{dt}t\lesssim||f||^r_{V_\theta}.$$
For this, we recall that $\f$ is the inverse of the function $\psi(x)=x^{\frac1p-\frac1q}\ummL x{-\a}p{\a}q,\quad x\in[\,0,1\,]$ so we have the following relations:
$$t=\f(t)^{\frac1p-\frac1q}\ummL{\f(t)}{-\a} p{\a}q,$$
$$1-\Log\f(t)\approx1-\Log t,$$
and
$$\dfrac{\psi'(x)}{\psi(x)}\approx\dfrac1x.$$
Therefore, we rewrite $II_2$ as:
$$II_2\approx\int_0^1\left[\f(t)^{(1-\theta)(\frac1p-\frac1q)}\ummL {\f(t)}{-\a(1-\theta)}p{\a(1-\theta)}q \sup_{\f(t)<s<1}\umL s{-\a}q \left(\int_{s}^1f_*^q(\tau)d\tau\right)^{\frac1q}\right]^r\dfrac{dt}t.$$
Making the change of variables, $t=\psi(x)\Longleftrightarrow x=\f(t)$
\begin{equation}\label{eqx57}
II_2\approx\int_0^1\left[x^{(1-\theta)(\frac1p-\frac1q)}\ummL x{-\a(1-\theta)}p{\a(1-\theta)}q\sup_{x<t<1}\umL  t{-\a}q\left(\int_t^1f_*^q(s)ds\right)^{\frac1q}\right]^r\dfrac{dx}x.
\end{equation}
Applying  Proposition~\ref{pxxx}, we have :
\begin{eqnarray*} & \int_0^1\left[x^{(1-\theta)(\frac1p-\frac1q)}\ummL x{-\a(1-\theta)}p{\a(1-\theta)}q \sup_{x<t<1}\umL t{-\a}q\left(\int_t^1f_*^q(s)ds\right)^{\frac1q}\right]^r\dfrac{dx}x\\
& \hskip1cm\LEQ c\int_0^1 \left[x^{\frac1{p_\theta}}\umL x{-\a}{p_\theta}f_*(x)\right]^r\dfrac{dx}x
\equiv c||f||^r_{V_\theta}.\end{eqnarray*}
\ \HF

From the preceding results, we can characterize the interpolation  spaces for small Lebesgue spaces.
\begin{theo}\label{tx}\ \\
Let $0<\theta<1,\ 1<r<+\infty,\ 1<p<q,\ \a>0$. Then
$$\Big(L^{(p,\a},L^{(q,\a}\Big)_{\theta,r}=L^{p_\theta,r}(\Log L)^{\frac\a{p'_\theta}}$$
with $\DST \dfrac1{p_\theta}=\dfrac{1-\theta}p+\frac\theta q,\ \dfrac1{p_\theta}+\dfrac1{p'_\theta}=1.$
\end{theo}
\pr
By the duality result on interpolation spaces, (see \cite{BK, Tartar}), we have, 
\begin{equation}\label{eqx70}
\Big[\Big(L^{(p,\a},L^{(q,\a}\Big)_{\theta,r}\Big]'=\Big(L^{p'),\a},L^{q'),\a}\Big)_{\theta,r'}=\Big(L^{q'),\a},L^{p'),\a}\Big)_{1-\theta,r'}
\end{equation} 
and noticing $\Big(L^{p'),\a}\Big)'=L^{(p,\a}$,
 $L^{\infty}(\O)$ is in $L^{(p,\a}\cap L^{(q,\a}$ and dense in each of these spaces.\\
From Theorem \ref{t2}, since $q'<p'$, we have\\
\begin{equation}\label{eqx71}
\Big(L^{q'),\a},L^{p'),\a}\Big)_{1-\theta,r'}=L^{q'_{1-\theta},r'}(\Log L)^{-\frac\a{q'_{1-\theta}}},
\end{equation}
with $\dfrac1{q'_{1-\theta}}=\dfrac\theta{q'}+\dfrac{1-\theta}{p'}=1-\dfrac1{p_\theta}=\dfrac1{p'_\theta}$.\\
Thus taking the associate space in the above equation (\ref{eqx71}) gives, taking into account (\ref{eqx70}),
$$\Big(L^{(p,\a},L^{(q,\a}\Big)_{\theta,r}=\Big[(L^{(p,\a},L^{(q,\a})_{\t,r}\Big]^{\prime\prime}=L^{p_\theta,r}(\Log L)^{\frac\a{p'_\theta}}$$
with $\dfrac1{p_\theta}=\dfrac{1-\theta}p+\dfrac\theta q,\ \dfrac1{p_\theta}+\dfrac1{p'_\theta}=1.$
\HF

A main consequence of Theorem \ref{t2} and Theorem \ref{tx} is the 
\begin{theo}\label{t3345}\ \\
Let $1<a<+\infty,\ \b\in\R,\ \b≠0,\ 1< r< +\infty$. Then the Lorentz-Zygmund space $L^{a,r}(\Log L)^\b$ is an interpolation space in the sense of Peetre of two Grand Lebesgue spaces if $\b<0$ and of two small Lebesgue spaces if $\b>0$.\end{theo}

\begin{rem}\ \\
If $\b=0$, it is  already known that the Lorentz space $L^{a,r}(\O)$ is an interpolation space of two classical Lebesgue spaces.\end{rem}
\section{Small Lebesgue space as interpolation of usual Lebesgue spaces}
The next proposition has been already proved on \cite{CC}, we drop its proofs :
\begin{propo}\label{p4.1}\cite{CC}\ \\
Let $1<p<+\infty$, $\a>0$. Then
$$L^{(p,\a}=\Big(L^p,L^\infty\Big)_{0,1;-\frac\a p+\a-1}.$$
\end{propo}

\begin{propo}\label{p4.2}\ \\ Let $1<p<q$, $\a>0$. Then 
$$L^{(p,\a}=\Big(L^p,L^q\Big)_{0,1;-1+\a-\frac\a p}.$$
\end{propo}
\pr
It is similar to the above Proposition \ref{p4.1}. Indeed, let set $W_{pq}$ the space on the RHS, then for $f\in W_{pq}$
$$||f||_{W_{pq}}=\int_0^1(1-\Log t)^{-1-\frac\a p+\a}\K p q \dfrac{dt}t.$$
By the Holmstedt formula, we have
$$\K p q\approx\left(\int_0^{t^\s}f_*^p(s)ds\right)^{\frac1p}+t\left(\int_{t^\s}^1f^q_*(s)ds\right)^{\frac1q}\hbox{
with }\dfrac1\s=\dfrac{q-p}{pq}.$$
Then, we deduce the equivalent expression of the norm:
\begin{equation}\label{eqx80}
||f||_{W_{pq}}\approx||f||_{L^{(p,\a}}+\int_0^1(1-\Log t)^{-1-\frac\a p+\a}t^{\frac1\s}\left(\int_t^1f_*^q(s)ds\right)^{\frac1q}\dfrac{dt}t.
\end{equation}

Let us show that the last term in relation (\ref{eqx80}) is less or equal to the norm of $f$ in $L^{(p,\a}$. Let us temporarily set 
$\b=-1+\a-\dfrac\a p$, let $1<\eps<\dfrac q\s+1$, then
\begin{eqnarray}
O_1&=&\int_0^1\left[(1-\Log t)^{\b q}t^{\frac q\s-\eps+1}t^{\eps-1}\int_t^1f^q_*(x)dx\right]^{\frac1q}\dtt\NN\\
&\LEQ&\int_0^1\left[(1-\Log t)^{\b q}t^{\frac q\s-\eps+1}\int_t^1x^{\eps-1}f_*^q(x)dx\right]^\frac1q\dtt.\label {eqx7000}
\end{eqnarray}
Applying Hardy inequality (see Theorem \ref{BN} and \cite{BR} Theorem 6.4), we obtain from the relation (\ref{eqx7000})
\begin{equation}\label{eqx7001}
O_1\lesssim\int_0^1\Big[(1-\Log t)^\b t^{\frac1\s+\frac1q}f_*(t)\Big]\dtt.
\end{equation}
But $\dfrac1\s+\dfrac1q=\dfrac1p$ then $\DST t^{\frac1\s+\frac1q}f_*(t)\LEQ\left(\int_0^t f^p_*(x)dx\right)^{\frac1p}.$ So we derive
\begin{equation}\label{eqx82}
O_1\lesssim\int_0^1\left[(1-\Log t)^\b\left(\int_0^tf_*^p(x)dx\right)^{\frac1p}\right]\dtt=||f||_{L^{(p,\a}}.
\end{equation}
Thus, we have from relations (\ref{eqx80}) to (\ref{eqx82}):
$$||f||_{L^{(p,\a}}\lesssim||f||_{W_{pq}}\lesssim||f||_{L^{(p,\a}}.$$
\ \HF
\section{Interpolation between small, Grand Lebesgue spaces and the associated $K$-functional}
In this section, we want to determine the interpolation space   $\big(L^{(p,\a},L^{q),\b}\big)_{\theta,r}. $ Due 
to the technical aspect of the proof, we shall only consider  the case where $\a=\b=1$, the argument remains the same in the general case.\\
 We want to show the following theorem:
\begin{theo}\label{t5}\ \\
Let $0<\theta<1,\ 1\LEQ r<+\infty,\ p<q$. Then
$$\big(L^{(p},L^{q)}\big)_{\theta,r}=L^{p_\theta,r}\big(\Log L\big)^{\a_\theta}$$
where $\dfrac1{p_\t}=\dfrac{1-\t}p+\dfrac \t q,\ \ \a_\t=1-\theta-\dfrac1{p_\t}$.
\end{theo}
We shall need the following $K$-functional. 
\begin{theo}\label{t6}\ \\
Let $1<p<q<+\infty$. Then for all $t>0,\ f\in L^{(p}+L^{q)}$
\begin{eqnarray*}
\K {(p} {q)}&\approx&\int_0^{\f(t)}\umL s {-1} p\left(\int_0^sf_*^p(\tau)d\tau\right)^{\frac1p}\dfrac{ds}s
+\umL t{p-1}p\left(\int_0^{\f(t)}f_*^p(\tau)d\tau\right)^{\frac1p}\\
&&+t\sup_{\f(t)<s<1}\umL s{-1}q\left(\int_s^1f_*^q(\tau)d\tau\right)^{\frac1q}\\
&\dot=&\sum_{i=1}^3K_i(t)
\end{eqnarray*}
where $\f$  is an invertible  function from $[0,1]$ into itself satisfying the equivalence 
$$\f(t)^{\frac1p-\frac1q}(1-\Log \f(t))^{\frac{p-q+pq}{pq}}\approx t.$$
\end{theo}
\pr
We can apply Lemma  ~\ref{lemBR} for $A_0=L^{(p}$ and $A_1=L^{q)}$. 
Let $\Phi_0(t)=\left\|\chi_{[0,t]}\right \|_{L^{(p}}$ and $\Phi_1(t)=\left\|\chi_{[0,t]}\right \|_{L^{q}}$.
As we have  
\[ \Phi_0(t)\approx t^{\frac{1}{p}}(1-\Log t)^{\frac{1}{p'}}, \quad  \Phi_1(t)\approx t^{\frac{ 1}{q}}(1-\Log t)^{-\frac{1}{q}},\] 
we have that $\psi(t)=\frac{\Phi_0(t)}{\Phi_1(t)}$..
As  $p<q$ the  conditions \eqref{C0}  \eqref{C1} and \eqref{C2} easily follow by using  the Proposition~\ref{p2} and Proposition~\ref{p2.3}. and from Lemma~\ref{lemBR} we obtain the result.\HF

Next, we want to prove Theorem \ref{t5}. We need to show that the norm of $f$ in $W\dot=(L^{(p},L^{q)})_{\t ,r}$ is equivalent to its norm in $V=L^{p_\t,r}(\Log L)^{\a_\t}$.\\
For commodity, we set $$\l=\t\left(\frac1p-\frac1q\right),\ \a_\t=1-\t-\dfrac1{p_\t},\ a=\l-\t,\ \l_1=\dfrac1{p_\t}-\dfrac1q=(1-\t)\left(\dfrac1p-\dfrac1q\right).$$
According to Theorem \ref{t6}, the expression of the norm of $f$ in $W$ is composed with 3 terms:
$$||f||^r_W=\int_0^1\Big(t^{-\t}\K{(p}{q)}\Big)^r\dfrac{dt}t\approx N_1+N_2+N_3$$
with $N_i$ corresponding to the function $K_i$ and which reads as follow after a change of variables
\begin{eqnarray*}N_1&=&\int_0^1\left[t^{-\l}(1-\Log t)^a\int_0^t\umL s{-1}p\left(\int_0^sf_*^p(x)dx\right)^{\frac1p}\dfrac{ds}s\right]^r\dfrac{dt}t,\\
N_2&=&\int_0^1\left[t^{-\l}(1-\log t)^{\a_\t}\left(\int_0^tf_*^p(x)dx\right)^{\frac1p}\right]^r\dfrac{dt}t,\\
N_3&=&\int_0^1\left[t^{\l_1}(1-\Log t)^{\a_\t+\frac1q}\sup_{t\LEQ s<1}\umL s{-1}q
\left(\int_s^1f_*^q(x)dx\right)^{\frac1q}\right]^r\dfrac{dt}t.\end{eqnarray*}
We start with a lower bound for $||f||_W$.
\begin{lem}\label{l5.2.1}\ \\
One has
$$||f||_W^r\gtrsim N_3\gtrsim ||f||_V^r=\int_0^1\left[t^{\frac1{p_\t}}(1-\Log t)^{\a_\t}f_*(t)\right]^r\dfrac{dt}t.$$\end{lem}
\pr
One has
\begin{eqnarray*}
N_3&\GEQ&\int_0^1\left[t^{\l_1}(1-\Log t)^{\a_\t}\left(\int_t^1f_*^q(x)dx\right)^{\frac1q}\right]^r\dfrac{dt} t\\
&\GEQ&\int_0^{\frac12}\left[t^{\l_1}(1-\Log t)^{\a_\t}\left(\int_t^{2t}f_*^q(x)dx\right)^{\frac1q}\right]^r\dfrac{dt}t\\
&\gtrsim&\int_0^{\frac12}\left[t^{\l_1+\frac1q}(1-\Log t)^{\a_\theta}f_*(2t)\right]^r\dfrac{dt}t.
\end{eqnarray*}
Thus, we have, after a change of variables
$$\int_0^1\left[t^{\frac1{p_\t}}(1-\Log t)^{\a_\t}f_*(t)\right]^r\dfrac{dt}t\lesssim N_3\lesssim||f||_W^r.$$
\ \HF

For the upper bound, we start with the estimate of $N_2$
\begin{lem}\label{l5.3}\ \\
One has $$N_2\lesssim||f||^r_V.$$
\end{lem}
\pr
We set $b=\dfrac rp\in]\,0,+\infty[,\ \b=p\a_\t$, we  note that $0<\l p<1$ we can write
$$N_2=\int_0^1\left[t^{-\l p}(1-\Log t)^\b\int_0^tf_*^p(x)dx\right]^b\dfrac{dt}t.$$
We may apply the Hardy's inequalities  of Theorem \ref{t5.3}  to derive
\begin{equation}\label{eqx200}
N_2\lesssim \int_0^1\Big[t^{1-\l p}(1-\Log t)^\b f_*^p(t)\Big]^b\dfrac{dt}t.
\end{equation}
The right hand side (RHS) of this inequality is equal to
\begin{equation}\label{eqx200b}
\int_0^1\Big[t^{\frac1p-\l}(1-\Log t)^{\a_\t}f_*(t)\Big]^r\dfrac{dt}t
\end{equation}
since $\dfrac1p-\l=\dfrac{1-\t}p+\dfrac\t q=\dfrac1{p_\t}$, we get the result.\HF

Next, we want to estimate $N_1$
\begin{lem}\label{l5.5}\ \\
One has:
$$N_1\lesssim \int_0^1\Big[t^{\frac1{p_\t}}(1-\Log t)^{\a_\t-1}f_*(t)]^r\dfrac{dt}t\LEQ ||f||^r_V.$$
\end{lem}
\pr
Let us set  $$H(s)=s^{-1}(1-\Log s)^{-\frac1p}\left(\int_0^s f_*^pdx\right)^{\frac1p}.$$
Then, from the Hardy's inequality given in Theorem \ref{t5.3} and the expression of $N_1$:
\begin{equation}\label{eqx300}
N_1\lesssim\int_0^1\Big[t^{1-\l}(1-\Log t)^aH(t)\Big]^r\dfrac{dt}t.
\end{equation}
Replacing $H$, noticing that $a-\dfrac1p=-\t-\dfrac1{p_\t}=\a_\t-1$, we then have from relation (\ref{eqx300})
\begin{equation}\label{eqx301}
N_1\lesssim\int_0^1\left[t^{-\l}(1-\Log t)^{\a_\t-1}\left(\int_0^tf_*^p dx\right)^{\frac1p}\right]^r\dfrac{dt}t.
\end{equation}
Applying again the Hardy inequality as in Lemma \ref{l5.3} we have
\begin{equation}\label{eqx302}
N_1\lesssim\int_0^1\Big[t^{\frac1p-\l}(1-\Log t)^{\a_\t-1}f_*(t)\Big]^r\dfrac{dt}t.
\end{equation}
This gives the result.\HF

It remains to estimate the term $N_3$. We have

\begin{lem}\label{l5.6}\ \\
One has $$N_3\lesssim||f||_V^r.$$
\end{lem}
The key lemma to estimate $N_3$ is the analogous of Proposition \ref{pxxx}
\begin{propo}\label{p5.1}\ \\ 
Let $K_d$ be a decreasing nonnegative function on $(0,1)$, $\nu>0$ and $\b$ two real numbers. Then
$$\int_0^1\Big[t^\nu(1-\Log t)^\b\sup_{t\LEQ s<1}(1-\Log s)^{-\frac 1q}K_d(s)\Big]^r\dfrac{dt}t\approx\int_0^1\Big[t^\nu(1-\Log t)^{\b-\frac1q}K_d(t)\Big]^r\dfrac{dt}t.$$
\end{propo}
The proof follows the same argument as for Proposition \ref{pxxx}.\\
Applying the Proposition \ref{p5.1}, we deduce 
\begin{equation}\label{eqx400}
N_3\lesssim\int_0^1\left[t^{\l_1}(1-\Log t)^{\a_\t}\left(\int_t^1f_*^q(x)dx\right)^{\frac1q}\right]^r\dfrac{dt}t.
\end{equation}

(The function $\DST K_d(t)=\left(\int_t^1 f_*^q(x)dx\right)^{\frac1q}$ is decreasing).

Again applying Hardy inequalities (see Theorem \ref{BN}) we have
\begin{equation}\label{eqx401}
N_3\lesssim \int_0^1\Big[t^{\frac1q +\l_1}(1-\log t)^{\a_\t}f_*(t)\Big]^r\dfrac{dt}t.
\end{equation}
But  we have $\DST\frac1q+\l_1=\dfrac1{p_\t}$ so the RHS of the above equation (\ref{eqx401})  is the norm of $f$ in $V$ at the power $r$.\HF\\
We have shown
$$||f||_W^r\lesssim N_1+N_2+N_3\lesssim||f||^r_V$$ and $$||f||_W^r\gtrsim N_3\gtrsim ||f||^r_V.$$
This proves Theorem \ref{t5}.\HF

\section{The critical case $p=q$. The interpolation space $(L^{p)},L^{(p})
$ and its $K$-functional}

The preceding study can be extended to the case where $p=q$. In this case, the inverse function of $\psi_1(x)=(1-\Log x)^{-1},$  say $\f_1(t)=e^{1-\frac1t}$ will play the fundamental role to express the $K$-functional. Note that in this case we can't use the Lemma~\ref{lemBR}. we should do it it by direct calculation. . 

\begin{theo}\label{t6.1}{\bf $K$-functional for $(L^{p)},L^{(p})$}\\
For $1<p<+\infty,\ 0<t<1,\ f\in L^{p)}+L^{(p}$, one has
$$\K{p)}{(p}\approx\!\!\!\!\!\!\sup_{0<s<e^{1-\frac1t}}\!\!\!\umL s{-1}p\left(\int_s^{e^{1-\frac1t}}\!\!\!\!\!\!f_*^p(x)dx\right)^{\frac1p}
+t\int_{e^{1-\frac1t}}^1\!\!\!\!\!\!\umL s{-1}p\left(\int_{e^{1-\frac1t}}^s\!\!f_*^p(x)dx\right)^{\frac1p}\dfrac{ds}s.$$
\end{theo}
\pr We only sketch it since the arguments are similar to the proof of Theorem \ref{t3}.\\ We consider 
$f=g+h\in L^{p)}+L^{(p}$. Then for $t>0,\ f_*(t)\LEQ g_*\left(\dfrac t2\right)+h_*\left(\dfrac t2\right)$. We derive
$$K^0(t)=\sup_{0<s<\e}\umL s{-1}p\left(\int_s^\e f_*^p(x)dx\right)^{\frac1p}\LEQ I_5+I_6.$$
with

$$I_5=\sup_{0<s<\f_1(t)}\umL  s{-1}p\left(\int_s^{\f_1(t)}g_*^p\xd dx\right)^\up.$$
With a change of variables, we have
\begin{equation}\label{eqx5500}
I_5\lesssim\sup_{0<s<\f_1(t)}\umL {\frac s2}{-1}p\left(\int_{\frac s2}^{\f_1(t)/2}g_*^p(x)dx\right)^\up\lesssim||g||_{L^{p)}}.
\end{equation}
While for the term
$$I_6=\sup_{0<s<\f_1(t)}\umL s{-1}p\left(\int_s^{\f_1(t)}h_*^p\xd dx\right)^\up,$$
since $1-\Log \f_1(t)=\dfrac1t$ and $s\to\umL s{-1}p$ is increasing, we have
\begin{equation}\label{5501}
I_6\LEQ t^\up\left(\int_0^{\f_1(t)}h_*^p(x)dx\right)^\up=t\times t^{\up-1}\left(\int_0^{\f_1(t)}h_*^p(x)dx\right)^\up.
\end{equation}
We replace the quantity $t^{\frac1p-1}$ by making an  integration of the following integral:
\begin{equation}\label{5502}
\int_{\f_1(t)}^1\umL s{-1}p\dss=\dfrac1{1-\frac1p}\Big[t^{\frac1p-1}-1\Big].
\end{equation}
Therefore, we have from (\ref{5501}) to (\ref{5502})
\begin{equation}\label{5503}
I_6\lesssim t\int_{\f_1(t)}^1\umL s{-1}p\left(\int_0^sh^p_*(x)dx\right)^\up\dss+t||h||_{L^p},
\end{equation}
\begin{equation}\label{5504}
I_6\lesssim t\Big(||h||_{L^{(p}}+||h||_{L^p}\Big)\lesssim t||h||_{L^{(p}}.
\end{equation}
Thus, we have
\begin{equation}\label{5505}
K^0(t)\lesssim ||g||_{L^{p)}}+t||h||_{L^{(p}}.
\end{equation}

While for the second term in the expression of $K$, say
$$K^1(t)=t\int_{\f_1(t)}^1\umL s{-1}p\left(\int_{\f_1(t)}^s f_*^p(x)dx\right)^\up\dss.$$
One has
$$K^1(t)\LEQ I_7+I_8,$$
with $$I_7=t\int_{\f_1(t)}^1\umL s{-1}p\left(\int_{\f_1(t)}^sg_*^p\xd dx\right)^\up\dss,$$
and 
$$I_8=t\int_\e^1\umL s{-1}p\left(\int_{\f_1(t)}^s h_*^p\xd dx\right)^\up\dss.$$

To estimate $I_7$, we use the relation $\DST\dfrac1t=1-\Log\f_1(t)$ to derive
\begin{equation}\label{5506}
I_7\lesssim t\left(\int_{\f_1(t)/2}^1g_*^p(x)dx\right)^\up t^{\frac1p-1} \lesssim t\umL{\dfrac{\f_1(t)}2}{p-1}p\left(\int_{\f_1(t)/2}^1g_*^p(x)dx\right)^\up
\lesssim||g||_{L^{p)}}.
\end{equation}
While for the second term $I_8$, we have after a change of variables
\begin{equation}\label{5507}
I_8\lesssim t\int_0^1\umL x{-1}p\left(\int_0^xh_*^p(\tau)d\tau\right)^\up\dfrac{dx}x=t||h||_{L^{(p}}.
\end{equation}
Then, we deduce from (\ref{5506}) to (\ref{5507})

\begin{equation}\label{5508}
K^1(t)\lesssim||g||_{L^{p)}}+t||h||_{L^{(p}}.
\end{equation}

Therefore, we have for all $f=g+h$
$$K^0(t)+K^1(t)\lesssim ||g||_{L^{p)}}+t||h||_{L^{(p}}.$$
That is
\begin{equation}\label{5509}
K^0(t)+K^1(t)\lesssim \K{p)}{(p}.
\end{equation}
For the reverse, we use the same decomposition as before $f\in L^{p)}+L^{(p}$, say
$$g=\Big(|f|-f_*(\f_1(t))\Big)_+=\Big(|f|-f_*(\f_1(t))\Big)\chi_{\{|f|>f_*(\f_1(t))\}},\quad h=f-g,$$ so that $h_*+g_*=f_*$ and
$$g_*=\Big(f_*-f_*(\f_1(t))\Big)_+,\quad h_*=f_*(\f_1(t))\chi_{_{(0,\f_1(t))}}+f_*(s)\chi_{(\f_1(t),1)}.$$

With those expressions, we derive
\begin{equation}\label{5510}
||g||_{L^{p)}}\LEQ\sup_{0<s<\f_1(t)}\umL s 1 p\left(\int_s^{\f_1(t)}f_*^p(x)dx\right)^\up= K^0(t).
\end{equation}

While for the term in $h$,
\begin{eqnarray}
t||h||_{L^{(p}}&\LEQ& t\left(\int_0^{\f_1(t)}\umL s{-1}p s^\up\dss\right)f_*(\f_1(t))\nonumber\\
&&+t\int_{\f_1(t)}^1\umL s{-1}p\dss\f_1(t)^\up f_*(\f_1(t))\nonumber\\
&&+t\int_{\f_1(t)}^1\umL s{-1} p\left(\int_{\f_1(t)}^sf_*^p(x)dx\right)^\up\dss.\label{5511}
\end{eqnarray}

For the first integral in (\ref{5511}), one can use Proposition \ref{p2} to derive
\begin{equation}\label{5512}
\int_0^{\f_1(t)}s^{\frac1p-1}\umL s{-1}pds\lesssim \f_1(t)^\up\umL{\f_1(t)}{-1}p.
\end{equation}
So we obtain
\begin{eqnarray}
t\int_0^{\f_1(t)}\umL s{-1}ps^\up\dss f_*(\f_1(t))&\lesssim& t\umL{\f_1(t)}{-1}p\left(\int_{\f_1(t)/2}^{\f_1(t)}f_*^p(x)dx\right)^\up\nonumber\\
&\lesssim&t\sup_{0<s\LEQ\frac{\f_1(t)}2}\umL s{-1}p\left(\int_s^{\f_1(t)}f_*^p(x)dx\right)^\up\nonumber\\
&\lesssim&t K^0(t)\LEQ K^0(t).\label{5513}
\end{eqnarray}
For the second integral, we use relation (\ref{5502}) to derive
\begin{equation}\label{5514}
B=t\int_{\f_1(t)}^1\umL s{-1}p\dss\f_1(t)^\up f_*(\f_1(t))\lesssim t^\up\f_1(t)^\up f_*(\f_1(t)).
\end{equation}
Since $\DST\int_{e^{-\frac1t}}
^{\f_1(t)} 
dx
=\dfrac{e-1}e\f_1(t)$, we then have
\begin{equation}\label{5515}
\f_1(t)^\up f_*(\f_1(t))\lesssim\left(\int_{e^{-\frac1t}}^{\f_1(t)}f_*^p(x)dx\right)^\up.
\end{equation}
So that relations (\ref{5514}) and (\ref{5515}) imply
\begin{equation}\label{5516}
B=t\int_{\f_1(t)}^1\umL s{-1}p\dss\f_1(t)^\up f_*(\f_1(t))\lesssim(1-\Log e^{1-\frac1t})^{-\frac1p}\left(\int_{e^{-\frac1t}}^{\f_1(t)}f_*^p(x)dx\right)^\up.
\end{equation}
Therefore, we obtain
\begin{equation}\label{5517}
B\lesssim\sup_{0<s<\f_1(t)}\umL s{-1}p\left(\int_s^{\f_1(t)}f_*^p(x)dx\right)^\up=K^0(t).
\end{equation}

The last term in relation (\ref{5511}) is equal to $K^1(t)$.\\
Combining these last relations (\ref{5510}) to (\ref{5517}), we come to
$$||g||_{L^{p)}}+t||h||_{L^{(p}}\lesssim K^0(t)+K^1(t),$$
which implies that
$$\K{p)}{(p}\lesssim K^0(t)+K^1(t)\lesssim \K{p)}{(p}.$$
\ \HF
We can now identify the interpolation space $(L^{p)},L^{(p})_{\t,r}$ as we state in

\begin{theo}\label{t6.2}\ \\
Let $1<p<+\infty,\ 0<\t<1,\ 1\LEQ r<+\infty$. Then
$Z_{\t,r}\dot=(L^{p)},L^{(p})_{\t,r} $ has the following equivalent norm: \\
For $f\in Z_{\t,r},\ \b_\t=\t-\dfrac1p-\dfrac1r$
\begin{itemize}

\item if $\t<\dfrac1p$ then
$$||f||_{Z_{\theta,r}}\approx\left[\int_0^1\left[\umLtbt\left(\int_t^1\fep\right)^\up\right]^r\dtt\right]^\ur.$$

\item if $\t>\dfrac1p$ then $Z_{\theta,r}=G\Gamma(p,r,w), w(t)=t^{-1}(1-\Log t)^{\bt r}$ and 
$$||f||_{Z_{\theta,r}}\approx \left[\int_0^1\left[\umLtbt\left(\int_0^t\fep\right)^\up\right]^r\dtt\right]^{\frac1r}\qquad
.$$

\item if $\t=\dfrac1p$
$$||f||_{Z_{\theta,r}}\approx\left[\sum_{k=0}^{+\infty}
\left(\int_{{2^{1-2^{k+1}}}}^{2^{1-2^k}}\fep\right)^{\frac rp}\right]^\ur.$$
In particular 
$$(L^{p)},L^{(p})_{\frac1p,p}=L^p\qquad\hbox{and\qquad}
(L^{p)},L^{(p})_{\t,p}=L^{p,p}(\log L)^{\t-\up}.$$
\end{itemize}
\end{theo}
\pr
According to Theorem \ref{t6.1}, we have for $f\in Z_{\t,r}$,
$$||f||^r_{Z_{\t,r}}\approx\int_0^1\left[t^{-\t}\sup_{0<s<\f_1(t)}\umL s{-1}p\left(\int_s^{\f_1(t)}f_*^p(y)dy\right)^\up\right]^r\dfrac{dt}t$$$$
+\int_0^1\left[t^{1-\t}\int_{\f_1(t)}^1\umL s{-1}p\left(\int_{\f_1(t)}^sf_*^p(y)dy\right)^\up\dss\right]^r\dfrac{dt}t\dot=F_1+F_2.$$

Using the change of variables $x=\f_1(t)$ so that $t=(1-\log x)^{-1}$, we have
$$F_1=\int_0^1\left[(1-\Log x)^\t\sup_{0<s<x}\umL s{-1}p\left(\int_s^xf_*^p(y)dy\right)^\up\right]^r\dfrac{dx}{(1-\Log x)x},$$
$$F_2=\int_0^1\left[(1-\Log x)^{\t-1}\int_x^1\umL s{-1}p\left(\int_x^sf_*^p(y)dy\right)^\up\dss\right]^r\dfrac{dx}{(1-\Log x)x}.$$
Let us start with the upper bound of this norm.
\begin{propo}\label{p6.1}\ \\
One has
$$F_1+F_2\lesssim\int_0^1\left[(1-\Log x)^{\t-\frac1r-\frac1p}\left(\int_0^xf_*^p(s)ds\right)^{\frac1p}\right]^r\dfrac{dx}x,\quad\hbox{if }\t<1.$$
\end{propo}
\pr
For the term $F_1$, we have
$$\left(\int_s^xf^p_*(t)dt\right)^{\frac1p}\LEQ\left(\int_0^x f_*^p(t)dt\right)^{\frac1p}\hbox{ for }s<x,$$
so, one deduces
\begin{equation}\label{eqx800}
F_1\LEQ\int_0^1\left[(1-\Log x)^{\t-\frac1r}\left(\int_0^xf_*^p(s)ds\right)^{\frac1p}\sup_{0<s<x}\umL s{-1}p\right]^r\dfrac{dx}x,
\end{equation}
which gives the upper estimate since $\DST\sup_{0<s<x}\umL s{-1}p=\umL x{-1}p$.\\
For the term $F_2$, we have
\begin{equation}\label{eqx801}
F_2\LEQ\int_0^1\left[(1-\log x)^{\t-1-\frac1r}\int_x^1\umL s{-1}p\left(\int_0^sf_*^p(y)dy\right)^{\frac1p}\dss\right]^r\dfrac{dx}x.
\end{equation}
Since $\t<1$, one may apply the Hardy inequality Theorem \ref{tBR6.5} (see \cite{BR} Theorem 6.5) to derive
\begin{equation}\label{eqx802}
F_2\LEQ\int_0^1\left[(1-\Log x)^{\t-\frac1r-\frac1p}\left(\int_0^xf_*^p(y)dy\right)^{\frac1p}\right]^r\dfrac{dx}x.
\end{equation}
\ \HF

\begin{propo}\label{p6.2}\ \\
For $0<\t<1$, we have
$$F_1+F_2\lesssim\int_0^1\left[(1-\Log x)^{\t-\frac1r-\frac1p}\left(\int_x^1f_*^p(y)dy\right)^\up\right]^r\dfrac{dx}x.$$
\end{propo}
\pr
We start with $F_2$. Since $$\int_x^sf_*^p(y)dy\LEQ\int_x^1f^p_*(y)dy,$$
the expression of $F_2$ can be estimated as :
\begin{equation}\label{eqx950}
F_2\LEQ\int_0^1\left[(1-\Log x)^{\t-1-\frac1r}\left(\int_x^1f_*^p(y)dy\right)^\up\left(\int_x^1(1-\Log s)^{-\up}\dss\right)\right]^r\dfrac{dx}x.
\end{equation}
By integration as for relation (\ref{5502}), one has
\begin{equation}\label{eqx951}
 \int_x^1(1-\Log s)^{-\up}\dss=\dfrac p{p-1}\left[(1-\Log x)^{1-\up}-1\right].
 \end{equation}
Therefore, we deduce from (\ref{eqx950}) and  (\ref{eqx951}), the
\begin{equation}\label{eqx952}
 F_2\lesssim\int_0^1\left[(1-\Log x)^{\t-\frac1p-\frac1r}\left(\int_x^1f_*^p(y)dy\right)^\up\right]^r\dfrac{dx}x.
 \end{equation}
While for the second term $F_1$, we deduce from Bennett-Rudnick's Lemma ( \cite{BR} Lemma 6.1)
$$\sup_{0<s<x}(1-\Log s)^{-\up}\left(\int_s^1f^p_*(y)dy\right)^\up\LEQ 
c\int_0^x\left[(1-\Log s)^{-\up-1}\left(\int_s^1 f_*^p(y)dy\right)^\up\right]\dss.$$

Setting temporarily $\DST\psi(s)=(1-\Log s)^{-\up-1}\left(\int_s^1f_*^p(y)dy\right)^\up,$ we then deduce
$$F_1\lesssim\int_0^1\left[(1-\Log x)^{\t-\frac1r}\left(\int_0^x\psi(s)\dss\right)\right]^r\dfrac{dx}x.$$
By Hardy inequality Theorem \ref{tBR6.5} (see \cite{BR} Theorem 6.5).

$$F_1\lesssim\int_0^1\left[(1-\Log x)^{\t-\up-\frac1r}\left(\int_x^1f_*^p(y)dy\right)^\up\right]^r\dfrac{dx}x.$$
.\HF

For the lower bound, we need few lemmas

The first lemma is a consequence of a general lemma given in Goldman Heinig and Stepanov (\cite{GHS}, see also \cite{GP} Lemma 3.1).

\begin{lem}\label{l312016}\ \\
Let $t_k=2^{1-2^k}$, $\forall\,k\in\N$, $H$ a nonnegative locally integrable function on $(0,1)$,  $(q,\l)\in]0,+\infty[^2$. Then
\begin{enumerate}
\item $\DST\sum_{k\in\N}\left(
\int_0^{t_k}H(x)dx\right)^q2^{\l kq}
\approx\sum_{m\in\N}\left(2^{\l m}\int_{t_{m+1}}^{t_m}H(x)dx\right)^q,$

\item $\DST\sum_{k\in\N}\left(
\int_{t_{k+1}}^1H(x)dx\right)^q2^{-\l kq}\approx\sum_{m\in\N}\left(2^{-\l m}\int_{t_{m+1}}^{t_m}H(x)dx\right)^q.$
\end{enumerate}
\end{lem}
Here $\N=\Big\{0,1,2,\ldots\Big\}$ the set of natural numbers.\\
The next Lemma can be obtained by straightforward computation

\begin{lem}\label{l322016}\ \\
Let $t_k=2^{1-2^k},\ k\in\N,\ \l\neq0$. Then, one has
\begin{enumerate}
\item For any $s\in [t_{k+1},t_k]$ $$2^k\approx 1-\Log s,$$
\item $$\int_{t_{k+1}}^{t_k}(1-\Log t)^{\l-1}\dtt\approx2^{k\l},$$
for all $k\in\N$.
\end{enumerate}
\end{lem}

As a corollary of the above lemmas we have:

\begin{lem}\label{l332016}\ \\
Let $t_k=2^{1-2^k},\ k\in\N,\ \l\neq0, \ q>0$ and $H$ a nonnegative locally integrable function on $(0,1)$. Then
\begin{itemize}
\item If $\l>0$ then 
$$\sum_{k\in\N}\left(\int_0^{t_k}H(x)dx\right)^q2^{\l kq}\gtrsim\int_0^1\left[(1-\Log t)^\l\int_0^t H(x)dx\right]^q\dfrac{dt}{(1-\Log t)t}$$
and one has the equivalence if $\DST\int_0^1 H(x)dx\lesssim\int_0^{\frac12}H(x)dx$ (for instance $H$ decreasing).
\item If $\l<0$
$$\sum_{k\in\N}\left(\int_{t_{k+1}}^1H(x)dx\right)^q2^{\l kq}\approx\int_0^1\left[(1-\Log t)^{\l}\int_t^1H(x)dx\right]^q\dfrac{dt}{(1-\Log t)t}.$$
\end{itemize}
\end{lem}
\pr
If $\l>0$, we use statement 2.) of Lemma \ref{l322016}
\begin{equation}\label{3321}
2^{\l kq}\approx\int_{t_{k+1}}^{t_k}(1-\Log t)^{\l q-1}\dtt.
\end{equation}
Thus, we have
\begin{equation}\label{eq3322}
\sum_{k\in\N}\left(\int_0^{t_k} H(x)dx\right)^q 2^{\l kq}
\approx\sum_{k\in\N}\left(\int_0^{t_k}\intH\right)^q\int_{t_{k+1}}^{t_k}(1-\Log t)^{\l q-1}\dtt.
\end{equation}
Since for $t<t_k$, $\DST\int_0^tH(x)dx\LEQ\int_0^{t_k}H(x)dx$,  then we derive from relation (\ref{eq3322})
\begin{equation}\label{eq3323}
\sum_{k\in\N}\left(\int_0^{t_k}\intH\right)^q2^{\l kq}\gtrsim\sum_{k\in\N}
\int_{t_{k+1}}^{t_k}\left(\int_0^t\intH\right)^q(1-\Log t)^{\l q-1}\dtt,
\end{equation}
and we have for any $G$ nonnegative
\begin{equation}\label{eq3324}
\sum_{k\in\N}\int_{t_{k+1}}^{t_k}G(t)dt=\int_0^1G(t)dt.
\end{equation}
This ends the proof of the lower bound of the first statement.\\

For the upper bound, we change the index of summation and use relation (\ref{3321})
\begin{eqnarray}
\sum_{k=1}^{+\infty}\left(\int_0^{t_k}\intH\right)^q 2^{\l kq}
&\approx&\sum_{j\in\N}\left(\int_0^{t_{j+1}}\intH\right)^q2^{\l jq}\nonumber\\
&\lesssim&\sum_{j\in\N}\int_{t_{j+1}}^{t_j}\left(\int_0^t\intH\right)^q (1-\Log t)^{\l q-1}\dtt\nonumber\\
&\lesssim&\int_0^1(1-\Log t)^{\l q-1}\left(\int_0^t\intH\right)^q\dtt.\label{eq3335}
\end{eqnarray}
Moreover, if $\DST\int_0^1\intH\lesssim\int_0^{\frac12}\intH$ and since
\begin{eqnarray}
\left(\int_0^{\frac12}\intH\right)^q
&\lesssim&\int_{\frac12}^1\left(\int_0^t\intH\right)^q(1-\Log t)^{\l q-1}t^{-1}dt\NN\\
&\lesssim&\int_0^1\left[(1-\Log t)^\l\int_0^t\intH\right]^q\dfrac{dt}{(1-\Log t)t}.\label{eq3336}
\end{eqnarray}
Then from the two relations (\ref{eq3335}) and (\ref{eq3336}), we obtain the upper bound for $\l>0$.\\
Same argument holds for $\l<0$ for having the second statement.

The above Lemma \ref{l332016} holds for $\l=0$, noticing that $\DST\int_{t_{k+1}}^{t_k}\dfrac{dt}{(1-\Log t)t}\approx 1$. Then under the same conditions as before, we have

$$\sum_{k\in\N}\left(\int_0^{t_k}\intH\right)^q\approx
\int_0^1\left[\int_0^t\intH\right]^q\dfrac{dt}{(1-\Log t)t} 
$$
and
$$\sum_{k\in\N}\left(\int_{t_{k+1}}^1\intH\right)^q\approx
\int_0^1\left[\int_t^1\intH\right]^q\dfrac{dt}{(1-\Log t)t}.$$
\vspace{-1.5cm}\ \HF\\
To obtain the lower bound for $F_1$, we first show the

\begin{lem}\label{l342016}\ \\
One has
$$F_1\gtrsim\sum_{k\in\N}2^{k(\t-\up)r}\left(\int_{t_{k+1}}^{t_k}\fep\right)^{\frac rp}
\hbox{ with } t_k=2^{1-2^k},\ k\in\N.$$
\end{lem}
\pr
Let us set $\DST G(x)=\sup_{0<s<x}(1-\Log s)^{-\frac rp}\left(\int_s^xf_*^p(y)dy\right)^{\frac rp},\ q=\dfrac rp$ and write
$$F_1=\sum_{k\in\N}\int_{t_{k+1}}^{t_k}G(x)(1-\Log x)^{\t r-1}\dfrac{dx}x.$$
Then, we derive from Lemma \ref{l322016} and the definition of $G$

\begin{equation}\label{eq3325}
F_1\gtrsim\sum_{k\in\N}G(t_k)\int_{t_{k+1}}^{t_k}(1-\Log x)^{\t r-1}\dfrac{dx}x\approx\sum_{k\in\N}G(t_k)2^{k\t r}.
\end{equation}
But $\DST G(t_k)=\sup_{i\GEQ k}\sup_{t_{i+1}<s<t_i}(1-\Log s)^{-q}\left(\int_s^{t_k}f^p_*(y)dy\right)^q$
and Lemma \ref{l322016} implies that for $t_{i+1}<s<t_i $,\ \  $(1-\Log s)^{-q}\approx 2^{-iq}$, so that 
$$ G(t_k)\approx\sup_{i\GEQ k}2^{-iq}\sup_{t_{i+1}<s<t_i}\left(\int_s^{t_k}f_*^p(y)dy\right)^q,$$
that is
\begin{equation}\label{eq3326}
G(t_k)\approx\sup_{i\GEQ k}2^{-iq}\left(\int_{t_{i+1}}^{t_k}f_*^p(y)dy\right)^q
\end{equation}
and then
\begin{equation}\label{eq3327}
G(t_k)\approx\sup_{i\GEQ k}2^{-iq}\left(\int_{t_{i+1}}^{t_k}f_*^p(y)dy\right)^q\GEQ2^{-kq}\left(\int_{t_{k+1}}^{t_k}f_*^p(y)dy\right)^q.
\end{equation}
Combining relation (\ref{eq3327}) and (\ref{eq3325}), we deduce 
$$ F_1\gtrsim\sum_{k\in\N}2^{k\t r-k\frac rp}\left(\int_{t_{k+1}}^{t_k}f_*^p(y)dy\right)^{\frac rp}.$$
This ends the proof.
\HF
As a corollary of the above lemmas, one has
\begin{theo}\label{t622016}\ \\

On has for $f\in\Ztr$
\begin{enumerate}
\item If $\t>\dfrac1p$ then
$$||f||_\Ztr^r\gtrsim\int_0^1\left[(1-\Log t)^{\t-\frac1p}\left(\int_0^tf^p_*(x)dx\right)^\up\right]^r\dfrac{dt}{(1-\Log t)t}.$$
\item If $\t<\dfrac1p$ then
$$||f||_\Ztr^r\gtrsim\int_0^1\left[(1-\Log t)^{\t-\frac1p}\left(\int_t^1f_*^p(x)dx\right)^\up\right]^r\dfrac{dt}{(1-\Log t)t}.$$
\end{enumerate}
\end{theo}
\pr
Let us set $q=\dfrac rp$, $\l=\left|\t-\dfrac1p\right|p$.\\
If $\t-\frac1p>0$, we may apply statement 1) of Lemma \ref{l312016} with $H(x)=f_*^p(x)$ to derive from Lemma \ref{l342016}
\begin{equation}\label{eq3328}
F_1\gtrsim \sum_{k\in\N}\left(2^{\l k}\int_{t_{k+1}}^{t_k}f^p_*(y)dy\right)^q\approx\sum_{k\in\N}\left[
\int_0^{t_k}f^p_*(x)dx\right]^q2^{\l kq}.
\end{equation}
So that relation (\ref{eq3328}) gives with the help of Lemma \ref{l332016}
$$F_1\gtrsim\sum_{k\in\N}\left(\int_0^{t_k}f^p_*(x)dx\right)^q2^{\l kq}
\gtrsim\int_0^1\left[(1-\Log t)^\l\int_0^tf^p_*(x)dx\right]^q\dfrac{dt}{(1-\Log t)t}.$$
This last inequality implies the first statement of Theorem \ref{t622016}, noticing that $$||f||^r_{\Ztr}\gtrsim F_1.$$

While for the second statement, we apply the statement 2) of Lemma \ref{l312016} to derive from Lemma \ref{l342016}:
$$F_1\gtrsim\sum_{k\in\N}\left(2^{-\l k}\int_{t_{k+1}}^{t_k}f^p_*(y)dy\right)^q\approx\sum_{k\in\N}\left[
\int_{t_{k+1}}^1f^p_*(x)dx\right]^q2^{-\l kq}.$$
So that
$$F_1\gtrsim\sum_{k\in\N}\left(\int_{t_{k+1}}^1f_*^p(x)dx\right)^q2^{-\l kq}\gtrsim\int_0^1\left[(1-\Log t)^{-\l}\int_t^1f^p_*(x)dx\right]^q\dfrac{dt}{(1-\Log t)t}.$$
The last inequality comes from the second statement of Lemma \ref{l332016} and this ends the proof of Theorem \ref{t622016} since $||f||_\Ztr^r\gtrsim F_1$.\HF

It remains to investigate in the particular case $\t=\dfrac1p$. The lower bound for $||f||_\Ztr$ comes from Lemma \ref{l342016}. It is sufficient to show the
\begin{lem}\label{l652016}\ \\
For any $\t\in]0,1[$, we have
$$\max(F_1;F_2)\LEQ\sum_{k\in\N}2^{kr(\t-\frac1p)}\left(\int_{t_{k+1}}^{t_k}f_*^p(s)ds\right)^{\frac rp}.$$
\end{lem}
\pr
The estimates for $F_1$ and $F_2$ follow the same argument, nevertheless we detail both estimates for clarity reason.

Following the notation in the proof of Lemma \ref{l342016}, we have
$$F_1=\sum_{k\in\N}\int_{t_{k+1}}^{t_k}G(x)(1-\Log x)^{\t r-1}\dfrac{dx}x.$$
 Following the above argument as in Lemma \ref{l342016}, we have
\begin{equation}\label{eq3340}
F_1\lesssim\sum_{k\in\N}G(t_k)2^{k\t r}\lesssim\sum_{k\in\N}2^{k\t r}\sup_{i\GEQ k} 2^{-iq}\left(\int_{t_{i+1}}^{t_k}f^p_*(y)dy\right)^q.
\end{equation}
Writing
$$\int_{t_{i+1}}^{t_k}f_*^p(y)dy=\sum_{j=k}^i  2^j2^{-j}\int_{t_{j+1}}^{t_j}f^p_*(y)dy,$$
one deduces
$$\sup_{i\GEQ k}2^{-iq}\left(\int_{t_{i+1}}^{t_k}f_*^p(y)dy\right)^q
\LEQ \sup_{i\GEQ k}2^{-iq}
\left(\sum_{j=k}^i2^j\right)^q \sup_{j\GEQ k}2^{-jq}\left(\int_{t_{j+1}}^{t_j}f^p_*(x)dx\right)^q$$

Therefore, we have
$$F_1\lesssim \sum_{k\in\N}2^{k\t r}\sup_{j\GEQ k}2^{-jq}\left(\int_{t_{j+1}}^{t_j}f^p_*(x)dx\right)^q.$$
So, setting $\DST A_j\dot=\left(\int_{t_{j+1}}^{t_j}f_*^p(y)dy\right)^q$ estimating the supremum term by a sum, one has 
$$
F_1\lesssim\sum_{k\in\N}2^{k\t r}\sum_{j=k}^{+\infty}2^{-jq} A_j=\sum_{j\in\N}2^{-jq}A_j\sum_{k=0}^j 2^{k\t r}.
$$
$$
F_1\lesssim\sum_{j\in\N}2^{j\t r-jq}A_j=\sum_{j\in\N}2^{j(\t-\frac1p)r}\left(\int_{t_{j+1}}^{t_j}f_*^p(y)dy\right)^q.
$$ 
For estimating $F_2$, we set $G_0(x)=\DST\left[\int_x^1\umL s {-1} p
\left(\int_x^sf_*^p(y)dy\right)^\up\dss\right]^r$ and write
$$F_2=\sum_{k\in\N}\int_{t_{k+1}}^{t_k}(1-\Log x)^{(\t-1)r-1}G_0(x)\dfrac{dx}x.$$
Then,  using Lemma \ref{l322016} and the fact that $G_0(x)\LEQ G_0(t_{k+1})$ for $t_{k+1}<x<t_k$, one has
\begin{equation}\label{eq5000}
F_2\lesssim\sum_{k\in\N}G_0(t_{k+1}) 2^{kr(\t-1)}.
\end{equation}
Writing $G_0(t_{k+1})$ as
$$G_0(t_{k+1})=\left[\sum_{i=0}^k\int_{t_{i+1}}^{t_i}\umL s{-1}p\left(\int_{t_{k+1}}^sf_*^p(y)dy\right)^\up\dss\right]^r,$$
and using Lemma \ref{l322016}, we then have
\begin{equation}\label{eq5001}
G_0(t_{k+1})\lesssim\left[\sum_{i=0}^{k}\left(\int_{t_{k+1}}^{t_i}f_*^p(y)dy\right)^\up2^{(1-\up)i}\right]^r.
\end{equation}
Using Lemma \ref{l312016} with an adequate step function $H$, 
the RHS of relation (\ref{eq5001})   can be estimated as
\begin{eqnarray*}
\sum_{i=0}^k\left(\int_{t_{k+1}}^{t_i}f_*^p(y)dy\right)^\up2^{(1-\up)i}
&=&\sum_{i=0}^k2^{(1-\up)i}\left(\sum_{j=i}^k\int_{t_{j+1}}^{t_j}f_*^p(y)dy\right)^\up\\
&\approx&\sum_{i=0}^k 2^{(1-\up)i}\left(\int_{t_{i+1}}^{t_i}f_*^p(y)dy\right)^\up.
\end{eqnarray*}
Summarizing the above relations, we come to
$$F_2\lesssim\sum_{k\in\N}2^{kr(\t-1)}\left[\sum_{i=0}^k2^{(1-\up)i}\left(\int_{t_{i+1}}^{t_i}f_*^p(y)dy\right)^\up\right]^r,$$
using the second statement of Lemma \ref{l312016}, by taking $H$ a step function on $[0,1]$ such that $\DST\int_{t_{i+1}}^{t_i} H(x)dx=2^{(1-\up)i}\left(\int_{t_{i+1}}^{t_i}f_*^p(y)dy\right)^\up , $ we deduce\ from this last inequality\
$$F_2\lesssim\sum_{i\in\N}2^{i(\t-\up)r}\left(\int_{t_{i+1}}^{t_i}f_*^p(y)dy\right)^{\frac rp}.$$
This ends the proof.\HF
\vspace{-1cm}

\begin{coro}\label{c1l6534}
{\bf of Lemma \ref{l652016} and Lemma \ref{l342016}}\\
One has for any $\t\in]\,0,1\,[$ and $f\in \Ztr$
$\DST ||f||^r_\Ztr\approx\sum_{k\in\N}2^{kr(\t-\frac1p)}\left(\int_{t_{k+1}}^{t_k}f^p_*(y)dy\right)^{\frac r p}.$
\end{coro}
As a consequence of the above corollary, we have
\begin{theo}\label{t45416}\ \\
One has
$$||f||_\Ztr\approx\left[\int_0^1(1-\Log t)^{\t r}\left(\int_0^t(1-\Log x)^{-1}f_*^p(x)dx\right)^{\frac rp}\dfrac{dt}{(1-\Log t)t}\right]^{\frac1r}.$$
\end{theo}
\pr
From the Lemma \ref{l322016}
$$\int_{t_{k+1}}^{t_k}2^{-k}f_*^p(y)dy\approx\int_{t_{k+1}}^{t_k}(1-\Log y)^{-1}f_*^p(y)dy,$$
and from the above Corollary \ref{c1l6534}
\begin{eqnarray*}
||f||^r_\Ztr&\approx&\sum_{k\in\N}2^{kr\t}\left(\int_{t_{k+1}}^{t_k}(1-\Log x)^{-1}f_*^p(x)dx\right)^{\frac rp}\\
&\approx&\sum_{k\in\N}2^{kr\t}\left(\int_0^{t_k}
(1-\Log x)^{-1}f_*^p(x)dx\right)^{\frac rp}\hbox{ (using Lemma \ref{l312016})}\\
&\approx&\int_0^1(1-\Log x)^{r\t}\left(\int_0^t(1-\Log t)^{-1}f_*^p(x)dx\right)^{\frac rp}\dfrac{dt}{(1-\Log t)t}\hbox{ (using Lemma \ref{l332016})}.
\end{eqnarray*}
\vspace{-1cm}

\ \HF

{\bf Acknowledgement}  : This paper has been written during the visit of the third  author (A. Gogatishvili) at 
the University of Poitiers in April 2016; he wishes to thank all the members of the department  for their kind hospitality. \\

The research of A. Gogatishvili and T.Kopaliani was in part supported by 
the Grant no DI/9/5-100/13 of the Shota Rustaveli National Science 
Foundation, and grant no~217282, Operators of Fourier analysis in 
some classical and new function spaces. \\
The research of A.Gogatishvili 
was partially supported by the grant P201/13/14743S of the Grant agency 
of the Czech Republic and RVO: 67985840.\\

The second author has been partially supported by the Gruppo Nazionale  
per l’Analisi Matematica, la Probabilità e le loro Applicazioni  
(GNAMPA) of the Istituto Nazionale di Alta Matematica (INdAM), by  
Project Legge 5/2007 Regione Campania “Spazi pesati ed applicazioni al  
calcolo delle variazioni” and by Università degli Studi di Napoli  
Parthenope through the project “sostegno alla Ricerca individuale”  
(triennio 2015-2017).
\\

\end{document}